\documentclass[12pt]{amsart}
\usepackage{caption}

\usepackage[dvipsnames]{xcolor}
\usepackage{latexsym}
\usepackage{tikz, pgfplots}
\usetikzlibrary{cd}
\usepackage{hyperref}
\theoremstyle{definition}

\usepackage{stmaryrd}

\usepackage{tkz-fct}  
\definecolor{bl}{RGB}{47,49,49}
\definecolor{skBlue}{RGB}{138 186 211}
\definecolor{egg}{RGB}{248,241,229}
\definecolor{citr}{RGB}{249,186,50}

\usepackage[papersize={210mm, 297mm}]{geometry}

\renewenvironment{frame}{}{}

\newcommand{\befr}{\begin{frame}}
\newcommand{\efr}{\end{frame}}
\usepackage{amsfonts,amssymb} 
\usepackage{stmaryrd}
\usepackage{amsmath}
\usepackage{mathtools}

\usepackage{cancel} 

\newtheorem{theorem}{Theorem}[section]

\newtheorem{proposition}[theorem]{Proposition}
\newtheorem{remark}[theorem]{Remark}
\theoremstyle{definition}
\newtheorem{definition}[theorem]{Definition}

\theoremstyle{remark}

\usepackage{xcolor}
\usepackage{xparse}
\usepackage{tikz}
\usetikzlibrary{matrix}
\usetikzlibrary{cd}

\definecolor{amber}{rgb}{1.0, 0.49, 0.0}
\definecolor{amethyst}{rgb}{0.45, 0.31, 0.59}
\definecolor{applegreen}{rgb}{0.55, 0.71, 0.0}
\definecolor{asparagus}{rgb}{0.0, 0.42, 0.24}
\definecolor{darkbyzantium}{rgb}{0.36, 0.22, 0.33}
\definecolor{darklavender}{rgb}{0.45, 0.31, 0.59}

\allowdisplaybreaks

\begin{document}
\title[higher Bessel kernels]
{Product formulas for the Higher Bessel functions}

\author{Ilia Gaiur, Vladimir Rubtsov and Duco van Straten}

\maketitle

 \begin{flushright}{\em To Maxim Kontsevich with our admiration.}
 \end{flushright}
\begin{abstract}
  We consider the generating function $\Phi^{(N)}$ for the reciprocals $N$-th power of factorials.
  We show a connection of product formulas for such series with the periods for certain families
  of algebraic hypersurfaces. For these families we describe their singular loci. We show that these
  singular loci are given by zeros of the Buchstaber-Rees polynomials, which define $N$-valued group laws.
  We describe a generalized Frobenius method and use it to obtain special expansions for multiplication kernels
  in the sense of Kontsevich. Using these expansions we provide some experimental results that connect
  $N$-Bessel kernels and the hierarchies of the palindromic unimodal polynomials. We study the properties of
  such polynomials and conjecture positivity of their roots. We also discuss the connection with Kloosterman
  motives as a version of the mirror duality.
\end{abstract}

\section*{Introduction}
In this paper we study certain functions that appear in the study of the quantum cohomology for $\mathbb{P}^{N-1}$. These objects play a central role in the standard examples of the mirror symmetry for Fano varieties as formulated by Givental and Dubrovin, linking quantum cohomology on the symplectic side and Landau-Ginzburg models on the algebro-geometric side. Inspired by mirror symmetry Givental discovered that a certain generating function for Gromov-Witten invariants of the Fano variety can be expressed as a period function for a Landau-Ginzburg model. We refer to \cite{Dubrovin1996,dubrovin1999painleve,givental1997tutorial,givental1995homological,givental1996equivariant,iritani2011quantum, kasprzyk2022laurent} for an introductory reading, to name a few. In this paper, we establish another connection with period functions for certain Landau-Ginzburg type models. This connection arises from consideration of the product formulas for the $I$-functions of $\mathbb{P}^{N-1}$, which we refer to as $N$-Bessel functions.

Our plan is to exhibit the motivic nature of multiplication kernels on the very explicit example of the kernels for these $N$-Bessel functions, which have a Taylor series of the form
\begin{equation}
    \Phi_N (x) = \sum\limits_{i=0}^\infty \frac{x^i}{i!^N},
\end{equation}
and represent the unique analytic solution for the following differential equation:
\begin{equation}
    (\theta_x^N-x)\psi=0,\quad \theta_x=x\frac{d}{dx},
\end{equation}
which we call the $N$-Bessel equation. For $N=1$ we are dealing with the ordinary exponential, whereas for
$N=2$, this series may be expressed in terms of the modified Bessel functions as follows:
\begin{equation}
    \Phi_2(x) = I_0(2\sqrt{x})=J_0(2\sqrt{-x}).
\end{equation}
The $N$-Bessel equation and its holomorphic solution play a prominent role in mirror symmetry for
${\mathbb P}^{N-1}$ as it arises naturally from its quantum cohomology, defined in terms of
Gromov-Witten invariants. The relation $p^N=0$ of the ordinary cohomology ring
$H^*(\mathbb{P}^{N-1})=\mathbb{C}[p]/p^N$
is thereby deformed to the relation $p^N=q$ in
$QH^*({\mathbb P}^{N-1})$. The Dubrovin-Givental' connection \cite{Dubrovin1996,dubrovin1999painleve,givental1995homological}
translates this (via the Schr\"odinger prescription) $p\leftrightarrow \theta_x, \,q \leftrightarrow x\cdot$ into the $N$-Bessel equation
$\theta_x^N \psi=x\psi$, whose solutions can be represented as exponential (or oscillatory) thimble integrals
of the Laurent polynomial (Landau-Ginzburg potential)
\[ S_N:=X_1+X_2+\ldots +X_n+\frac{x}{X_1X_2\ldots X_N}.\]
In particular
\[\Phi_N(x)= \frac{1}{(2\pi i)^N}\oint\limits_T e^{S_N(X_1,X_2,\ldots,X_N;x)} \frac{dX_1}{X_1}\frac{dX_2}{X_2}\ldots \frac{dX_N}{X_N} \]
 where $T$ stands for the torus defined by $|X_i|=\epsilon_i$. In general, the thimble integrals can be expressed as
 Laplace integrals of periods integrals of the Laurent polynomial $S_N$, which is the reason
 that $S_N$ is often referred to as mirror dual to $\mathbb {P}^{N-1}$. For more information, we refer to
\cite{givental1996equivariant,iritani2011quantum, kasprzyk2022laurent}. 
For $N=2$ we are dealing with $\mathbb { P}^1$ and integrals of $e^{X+x/X}$. The corresponding finite field
sums over $X$ are none other than the classical Kloosterman sums \cite{katz1988gauss,xu2022bessel}.\\

In this paper we study  {\em product formulae} for the higher Bessel function of the type
\begin{equation}
    \Phi_N(x)\Phi_N(y) = \frac{1}{2\pi i}\oint\limits_{T}K_N(x,y|z) \Phi_N(z) \frac{dz}{z},
\end{equation}
where $T$ is a circle centered at $z=0$ with a sufficiently small radius and $K_N(x,y|z)$ is a so-called {\em multiplication kernel} \cite{kontsevich2021multiplication}. Moreover, in this paper we also look at more general product formulae
\begin{equation}
    \Phi_N(x_1)\Phi_N(x_2)\dots \Phi_N(x_m) = \frac{1}{2\pi i}\oint\limits_{T}K_N^{(m)}(x_1,x_2\dots, x_m|x_0) \Phi_N(x_0) \frac{dx_0}{x_0}.
\end{equation}

Recently, the concept of multiplication kernels and its place in the field of differential
equations was revisited in the papers \cite{golyshev2021non, kontsevich2021multiplication}.
The suggestion that such multiplication kernels are of motivic nature was the starting point
of our interest in the problem. In this paper we will exhibit the motivic nature of the
kernels $K_N^{(m)}$, more precisely describe them in terms of period functions of a specific
family of hypersurfaces over a parameter space $\mathbb{P}^m$, in line with Kontsevich's ideas.
Following the philosophy of  \cite{golyshev2021non}, one can consider the existence of such product
identities for some solutions of differential equations as a non-Abelian version of the classical
Abel's theorem, which expresses a sum of $g+1$ integrals as sum of $g$ integrals at algebraically
related arguments. While Abel's theorem in a general sense provides a $g+1$ to $g$ commutative structure
for the line bundles over genus $g$ curves (solutions of the first order differential equations),
the non-Abelian version leads to the integral identities on the space of flat sections for vector
bundles (solutions of the $N$-th order differential equations). Such non-Abelian integral identities
may be seen as a lift of the corresponding structures on the spectral curves of the corresponding
connection. This was also explained in \cite{kontsevich2021multiplication} by introducing the notion
of semi-classical kernel and its quantization.

Although such motivic multiplication kernels are probably most useful for the study of non-rigid differential
equations in a sense of N. Katz \cite{katz1996rigid}, already the in the $N$-Bessel case and more generally
rigid hypergeometric systems, this lead to some non-trivial results. Rather then dealing with genus $g$ Abelian
law lifts, we show how here how the $N$-Bessel kernels lift the Buchstaber-Rees $N$-valued group laws for
$\mathbb{C}/\mathfrak{S}_N$. 

\vspace{0.2cm}
Let us briefly describe the structure of the paper and give an exposition of results. This paper may be divided in two parts - first one contains theorems and theoretical results, while the second one is concentrated mostly on the connections to the other fields, computational and conjectural results. The first part consists of:
\begin{enumerate}
    \item In section \ref{sec:2Bes} we expose the already known results on $2$-Bessel function, corresponding product formula and its appearance in different fields of mathematics and theoretical physics.
    \item Section \ref{sec:LG1} is dedicated to the product of two $N$-Bessel functions. In particular, we show that corresponding multiplication kernel coincides with certain period of special Landau-Ginzburg models for any $N\in \mathbb{Z}_{\geq 2}$. Moreover, the singularity loci for these models is studied and given in closed algebraic form. 
    \item Section \ref{sec:LG2} extends the results of section \ref{sec:LG1} to the case of multi-product formulas. Corresponding Landau-Ginzburg models and singular loci are provided.
\end{enumerate}
The second part contains:
\begin{enumerate}
    \setcounter{enumi}{3}
     \item In section \ref{sec:per} we consider special $1$-parameter families of LG potentials which periods provide kernels for $N$-Bessel duplication formulas. We compute corresponding Picard-Fuchs operators, mirror maps and Yukawa couplings for certain values of $N$. We also study the local monodromy at $0$ for the corresponding Picard-Fuchs operators and propose a conjecture for its Jordan form.
     \item In section \ref{sec:def} we perform a special deformation of the duplication kernel and study the corresponding Picard-Fuchs equation for the $4$-Bessel function. We provide a deformed version of mirror map as well as a deformation of Yukawa coupling. 
     \item Section \ref{sec:Buch} gives a brief introduction into $N$-valued groups and special construction of $N$-valued groups by Buchstaber and Rees. We show that singular loci of the LG models for the Bessel kernels coincide with the so-called Buchstaber-Rees polynomials. \item In section \ref{sec:FrobNum} we introduce another approach for the kernel construction as a special power series. Using this construction we obtain new families of integer coefficient polynomials associated with $N$-Bessel kernels. We study properties of these polynomials and state a conjecture about their roots.
     \item In the last section we provide a computational results which allows consider kernels as solutions of the Cauchy problem for a $2$-parameter family of Picard-Fuchs operators in one variable.
\end{enumerate}

The second part of this paper is mostly inspired by classical approaches in computational algebraic geometry. One of our initial aims was to find the Calabi-Yau differential operators introduced in \cite{van2017calabi}, which appear as annihilators of the kernel functions restricted to some special one-dimensional subvarieties. Using standard approaches for Laurent polynomials (see \cite{coates2016quantum, akhtar2012minkowski}), we study corresponding Picard-Fuchs equations. Moreover, it seems that the results we obtain generalize some densities of random walks on the plane \cite{borwein2012densities}. This may be key to the nice integrality properties that hold for the Frobenius expansions of the Bessel kernels introduced in Section \ref{sec:FrobNum}.

{\bf Acknowledgements.} 
The authors are grateful to V. Buchstaber, S. Galkin, V. Golyshev, D. Gurevich, J. Fres\'an, B. Khesin, M. Kontsevich, D. Maslov, A. Mellit, A. Odesskii, C. Sabbah, V. Schechtman, P. Vanhove, A. Veselov, W. Zudilin for fruitful discussions and advice. Authors also greatly acknowledge discussions with members of International Group de Travail on Differential Equations in Paris. Our special thanks are addressed  to Wadim Zudilin for his proofread of the text draft. I.G. was supported by Cecil King Scholarship award of London Mathematical Society and IH\'ES.  The paper was started and developed  during IG and VR visits in IHES in the period of 2022-23, talks by DvS at the GdT and finished during VR visit in 2024. VR, IG and DvS greatly acknowledge IH\'ES  stimulating atmosphere and excellent working conditions. VR is thankful also to  the Centre Henri Lebesgue, Programme ANR-11-LABX-0020-0.
DvS also acknowledges the support by the Deutsche Forschungsgemeinschaft (DFG, German Research Foundation) through the Collaborative Research Center TRR 326 \textit{Geometry and Arithmetic of UniformizedStructures}, project number 444845124.

\section{$2$-Bessel function multiplication kernels}\label{sec:2Bes}

In this section we review some known results pertaining to the $N=2$ case, i.e. the classical Bessel function.
The story of the Bessel multiplication kernel starts with the classical Clausen duplication identity 
\begin{equation*}
    \left(\sum\limits_{n=0}^{\infty}\frac{x^n}{n!^2}\right)^2 = \sum\limits_{n=0}^{\infty}\frac{x^n}{n!^2}\binom{2n}{n}. 
\end{equation*}
The right hand side may be seen as the Hadamard convolution of the initial series with the Taylor expansion for $1/\sqrt{1-4x}$. In the integral form it reads as
\begin{equation}\label{eq:clausen}
    \Phi_2(x)^2 = \frac{1}{2\pi i}\oint\limits_T \frac{1}{\sqrt{1-4x/z}}\Phi_2(z)\frac{dz}{z}.
\end{equation}
The formula above is written for the local expansion. In order to prove the equality "globally" we can use the fact that $\Phi(x)^2$ is a solution of the second symmetric power of the $2$-Bessel equation (for example see \cite{chudnovsky2006computer} for definition of $Sym^2$).  In other words, we have
\begin{equation*}
    (\theta_x^3-4x\,\theta_x-2x)\Phi_2(x)^2 =0.
\end{equation*}
\begin{proposition}
    The integral in the right hand side of \eqref{eq:clausen} is a solution of the equation
    \begin{equation*}
        (\theta_x^3-4x\,\theta_x-2x)\psi(x)=0.
    \end{equation*}
\end{proposition}
Below we give an elementary and transparent proof of this fact which avoids the local expansion of the integral, but uses analytical properties of the duplication kernel.
\begin{proof}
  This can be shown by Eulers method of differentiation under the integral sign. It was made into a powerful
  algorithm by Almkvist and Zeilberger \cite{almkvist1990method}. 
    Let us denote $K(x,t)=\frac{1}{\sqrt{1-4x/t}}$. Consider an integral over any closed contour which doesn't depend on $x$
    $$
    I(x) = \oint K(x,t) \,\phi(t)\,\,\frac{dt}{t}.
    $$
    It is easily checked that the kernel function $K(x,t)$ satisfies 
    \begin{equation*}
        \theta_x K+\theta_tK=0,\quad 2xK+(t-4x) \theta_tK=0.
    \end{equation*}
    Let us for simplicity write 
	$$
	f_x = \theta_x f,\quad f_t = \theta_t f,
	$$
 and compute $I_{xxx},I_x, xI$ 
 $$
	I_{xxx} = \oint K_{xxx} \cdot \phi\,\frac{dt}{t} =\oint - K_{ttt} \cdot \phi \,\frac{dt}{t} = \oint K\cdot   \phi_{ttt} \,\frac{dt}{t} = \oint K  \cdot (t \phi)_t \,\frac{dt}{t},
	$$
	$$
	I_x = \oint K_x \cdot \phi \, \frac{dt}{t} =\oint -K_t \cdot \phi\, \frac{dt}{t} = \oint K\cdot  \phi_t \,\frac{dt}{t},
	$$
	$$
	2x I = \oint 2x \cdot K \cdot \phi \,\frac{dt}{t} = \oint (4x-t)\cdot  K_t\cdot \phi \,\frac{dt}{t}= \oint K\cdot(  (t\phi)_t-4x\,\phi_t)\,\frac{dt}{t}
	$$
 Combining these integrals together one gets
 \begin{equation*}
     [\theta_x^3-4x\,\theta_x-2x]\,I = \oint K\cdot((t \phi)_t -4x \, \phi_t -(t\phi)_t+4x\,\phi_t)\,\frac{dt}{t} = 0 
 \end{equation*}
 Choosing $T$ as the contour of integration, one can see that both $I(x)$ and $\Phi(x)^2$ are solutions of the same differential equation. 
\end{proof}


The above formula has tight links with the Sonine--Gegenbauer integral formula \cite{gegenbauer1883,sonine1880recherches}. The kernel obtained above can be seen as a special case of the Sonine--Gegenbauer formula specified for $\Phi_2(x)$. For the Bessel functions of $0$-degree the Sonine--Gegenbauer formula reads as \cite{watson1966treatise}
    \begin{equation}\label{eq:SG}
    J_0(u)J_0(v) 
    =\frac{1}{2\pi} \int\limits_{u-v}^{u+v}\frac{J_0(w)wdw}{\sqrt{u^4+v^4+w^4-2u^2v^2-2u^2w^2-2v^2w^2}},\quad u<v\in \mathbb{R}.
\end{equation}
Here we use the standard notation for the Bessel function $J_\nu(x)$ of the $\nu$th-degree as the solution regular
at $0$ of the equation
\begin{equation}
    \left(x\frac{d}{dx}\right)^2 \psi+(x^2-\nu^2)\psi =0,
\end{equation}
with the local expansion of the form
\begin{equation}
    J_\nu(x) = \sum\limits_{m\geq 0} \frac{(-1)^m}{m!\Gamma(m+\nu+1)}\left(\frac{x}{2}\right)^{2m+\nu}.
\end{equation}
Furthermore, one can recognize  the Heron-configuration, which allows us to interpret the integral
in the right hand side of \eqref{eq:SG}. Indeed, Heron's formula expresses the area $S(u,v,w)$ of
the triangle with sides $u,v,w$ as
\begin{equation*}
    16 S^2 = (u+v+w)(u+v-w)(u-v+w)(u-v-w) = u^4+v^4+w^4-2u^2v^2-2u^2w^2-2v^2w^2.
\end{equation*}
The zero locus of this polynomial coincides with the integration limits in the Sonine--Gegenbauer formula, which in the real situation coincide with the boundaries in the moduli space of triangles.

The Sonine--Gegenbauer formula naturally transfers to the multiplication formula for two $\Phi_2$ functions. In particular, by performing a change of variables one gets
\begin{equation*}
    \Phi_2(x)\Phi_2(y) = \int\limits_{\sqrt{x}-\sqrt{y}}^{\sqrt{x}+\sqrt{y}} \frac{1}{\sqrt{x^2+y^2+z^2-2(xy+xz+yz)}}\Phi_2(z)dz,\quad y<x.
\end{equation*}
As it was shown in \cite{golyshev2021non}, the polynomial under the square root in the denominator,
\begin{equation}\label{eq:Kallen2}
x^2+y^2+z^2-2(xy+xz+yz),
\end{equation}
appears to be a special limit of the Kontsevich polynomial \cite{kontsevich2009notes}, which provides a multiplication kernel for a Heun equation with unipotent monodromies. The polynomial \eqref{eq:Kallen2} also known as the {\it K\"all\'en function} \cite{mizera2022landau} in the theoretical physics community and may be written as $t$-discriminant of the polynomial
\begin{equation*}
 t^2+(x+y+ z)t+xy+xz+yz.
\end{equation*}

More evidence that the multiplication kernels for classical Bessel functions are deeply connected with periods of algebraic families comes from the multiple integrals of the Bessel functions. Using the theory of Sturm-Liouville operators, one can show that the Bessel functions are "orthogonal" in the following sense:
\begin{equation}
    \int\limits_{0}^{\infty} J_0(\lambda x) J_0(\lambda y) \lambda\, d\lambda = \delta(x-y),
\end{equation}
where delta is the Dirac delta function. This implies the following expression for the multiplication kernel
\begin{equation}
    P(x_1,x_2,\dots x_m| x_0) = \int\limits_{0}^{\infty} J_0(\lambda x_1) J_0(\lambda x_2)\dots J_0(\lambda x_m)J_0(\lambda x_0) \lambda\, d\lambda.
\end{equation}
Such type of integrals were actual objects of study by Sonine \cite{sonine1880recherches} and Gegenbauer \cite{gegenbauer1883}. One can see the motivic nature of the kernel $P$ already for $m=3$: in that case one gets an elliptic integral of the first kind (see \cite{watson1966treatise}, 13.46, formula (9))
\begin{equation}
    \int\limits_{0}^{\infty}\prod\limits_{i=1}^4 J_0(a_nt)dt=\frac{1}{\pi^2}\left\{\begin{array}{l}
        \frac{1}{\Delta}K\left(\frac{\sqrt{a_1a_2a_3a_4}}{\Delta}\right), \\
         \frac{1}{\sqrt{a_1a_2a_3a_4}}K\left(\frac{\Delta}{\sqrt{a_1a_2a_3a_4}}\right),
    \end{array}\right.
\end{equation}
where $K(k)$ denotes complete elliptic integral of the first kind with modulus $k$, and $\Delta$ reads as
\begin{equation}
16\Delta^2=\prod\limits_{n=1}^4(a_1+a_2+a_3+a_4-2a_n).
\end{equation}
Obviously this function coincides with the period of a special family of elliptic curves
parametrized by ${\mathbb P}^3$. The extension of such formulas, which may be evaluated in terms of
elliptic integrals was studied in \cite{bailey2008elliptic}.

More general multiple $2$-Bessel kernels $K_2^{(m)}(x_0,x_1,\dots x_m|t)$ carry a more complicated geometric structure. One may connect these kernels with the following family of Landau-Ginzburg potentials
\begin{equation}
\left(X_0+X_1+\ldots+X_m\right)\left(\frac{x_0}{X_0}+\frac{x_1}{X_1}+\ldots+\frac{x_m}{X_m}\right)=1/t.
\end{equation}
The geometry of such families were studied in papers by Verrill \cite{verrill1996root,verrill2004sums}.
For the case $m=5$, in papers by Hulek and Verrill \cite{hulek2005modularity} and
Candelas c.s. \cite{candelas2020one,candelas2023mirror} the algebraic geometry of the corresponding
family of 3-dimensional Calabi-Yau varieties is studied.

These families also appear in the context of Feynman integrals of {\em banana graphs}. The two-loop {\em sunrise integral} studied by Weinzierl and M\"uller-Stach appears for $m=3$ family of elliptic curves \cite{muller2011second}. The three-loop banana graph corresponds to $m=4$ the family of K3 surfaces by Bloch, Kerr and Vanhove \cite{bloch2015feynman,bloch2016local}, and Klemm, D\"uhr \cite{bonisch2022feynman,klemm2020loop} for arbitrary $m$ (see also \cite{broadhurst2016feynman}).

In the next section we show how similar Landau-Ginzburg type periods appear as multiplication kernels for the higher Bessel functions.

\section{Powers of binomial coefficients and corresponding Landau-Ginzburg model}\label{sec:LG1}
In the theory of Bessel functions the spectral parameter enters the eigenfunction by re-scaling the argument. Indeed, let us consider a spectral problem of the form:
\begin{equation*}
     \left(x\frac{d}{dx}\right)^N \psi =  \lambda x \psi.
\end{equation*}
One may see that the analytic solution at zero equals to $\Phi_N(\lambda \cdot x)$. Using the approach from \cite{kontsevich2021multiplication} (see section 5.1),  one may see that a multiplication kernel for $\Phi_N$ then reads as a generating function for the $N$-powers of binomial coefficients, i.e.
\begin{equation}\label{eq:ker}
    K_N(x,y,z) = \sum\limits_{j,k} \binom{j+k}{k}^N\,\frac{x^jy^k}{z^{j+k}}.
\end{equation}
The kernel is a homogenious function, i.e.  $K(x,y,z)=K(x/z,y/z)$, which may be represented as follows. Define
\begin{equation}\label{eq:LGFam1}
  W(x,y,z) := z-\prod\limits_{j=1}^{N-1}(1+Y_{j}) \left(x+\frac{y}{Y_1Y_2\dots Y_{N-1}}\right) \in \mathbb{Z}[x,y,z,Y_1^{\pm},Y_2^{\pm},\ldots,Y_{N-1}^{\pm}]
\end{equation}
The locus $W(x,y,z)=0$ can be considered as family of hypersurfaces in the torus $(\mathbb{G}_m)^{N-1}$
parametrized by $(x:y:z) \in \mathbb{P}^2$, i.e. $\subset \mathbb{P}^2 \times (\mathbb{C}^*)^{N-1}$.
    \begin{theorem}\label{th:prod}
    The kernel $K_N(x,y,z)$  \eqref{eq:ker} has the following representation as $N-1$-fold integral
    \begin{equation}\label{eq:DublPeriod1}
 K_N(x,y,z) = \frac{1}{(2\pi i)^{N-1}}\oint\limits_{T} \frac{1}{W(x/z,y/z,1)}\prod\limits_{j=1}^{N-1}\frac{\operatorname{d} Y_j}{Y_j}.
\end{equation}
\end{theorem}

\begin{proof} Consider a polynomial
    \begin{equation}
        F = (\alpha +\beta Y_0)\prod\limits_{j=1}^{N-1}(1+Y_j).
    \end{equation}
     The $k$-th power of $F$ is a polynomial of the form
     \begin{equation*}
         F^k = \sum \binom{k}{j_0}\binom{k}{j_1}\binom{k}{j_2}\dots  Y_0^{j_0} Y_1^{j_1}Y_2^{j_2}\dots  \alpha^{k-j_0}\beta^{j_0}.
     \end{equation*}
     If we put
     \begin{equation*}
         Y_0=1/(Y_1Y_2\dots Y_{N-1}),
     \end{equation*}
     then the constant term $[F^k]_0$ is a polynomial in $\alpha$ and $\beta$ with the appropriate
     $N$-th powers of $\binom{k}{m}$ as coefficients. The standard movement rewrites the constant term
     series as a torus integral
     \begin{equation}
         \pi (\alpha, \beta) = \frac{1}{(2\pi i)^{N-1}}\oint \limits_T\prod\limits_{j=1}^{N-1}\frac{\operatorname{d} Y_j}{Y_j}\frac{1}{1-F} {\Big\vert}_{Y_0=\frac{1}{Y_1Y_2\dots Y_{N-1}}}.
     \end{equation}
      so we see that the formal expansion \eqref{eq:ker} coincides with the integral \eqref{eq:DublPeriod1}
   by setting  $\alpha=x/z$ and $\beta=y/z$ in the integral $\pi(\alpha,\beta)$.
\end{proof}

 As we have expressed the kernel $K_N(x,y,z)$ as a parametric integral over a rational function, it
 satisfies a system of differential equations on $\mathbb{P}^2$, the Picard-Fuchs or Gauss-Manin system.
 As such it defines a singular locus where the solutions branch. This locus coincides with the projection to
 $\mathbb{P}^2$ of the critical points for the family \eqref{eq:LGFam1}. For this we can give a
 nice explicit expression.
 
\begin{theorem}\label{th:sing1} The singularities in $\mathbb{P}^2$ of the Picard-Fuchs system for $K_N$
   are the coordinate triangle  $xyz=0$ and a degree $N$ projective rational curve $R_N$ 
defined by a symmetric polynomial 
\[\Delta_N(x,y,z)=x^N+y^N+z^N+\ldots,\]
given by the property 
\begin{equation}\label{eq:factorLoci}
    \Delta_N(u^N,v^N,w^N)=\prod_{\omega,\eta} (u+\omega v+\eta w),
\end{equation}
where in the product $\omega$ and $\eta$ run over the $N$-th roots of unity.
\end{theorem}
\begin{proof}
    Singularities of the Picard-Fuchs operator coincide with the singular fiber of the family
    \begin{equation}
        \mathcal{Y}_{[x:y:z]}: \quad x \prod_{i=1}^{N-1}(1+Y_i) + y \prod_{i=1}^{N-1}(1+Y_i^{-1}) = z.
    \end{equation}
   Taking derivatives one obtains the following set of equations defining the Jacobian ideal:
    \begin{equation}
        \prod\limits_{i\neq j}^{N-1}(1+Y_i)\left(x-\frac{y}{Y_j}\prod\limits_{i=1}^{N-1}Y_i^{-1}\right)=0,\quad j=1\dots N-1.
    \end{equation}
    1) {\em $Y_i\neq -1$}. Then it follows that possible singularities are the diagonals in $\mathbb{G}_{N-1}$:
    $$
    Y_1=Y_2=\dots =Y_{N-1} = T.
    $$
    Assuming this, we get a system of equations
    \begin{equation}
        \left\{\begin{array}{l}
            T^Nx=y  \\
            (1+T)^{N-1}(xT^{N-1} +y)=z
        \end{array}\right.
    \end{equation}
    which is equivalent to
    \begin{equation}
        \left\{\begin{array}{l}
            T^Nx=y  \\
            (1+T)^{N}y=T^Nz
        \end{array}\right.
    \end{equation}
    After the change of variables
    $$
    x=u^N,\quad y=v^N,\quad z=w^N,
    $$
    the system may be linearized in $T$ and reads as
    \begin{equation}
       \left\{ \begin{array}{l}
       Tu=v\cdot \chi\\
       (1+T)v=T \cdot w\cdot \eta \end{array}\right.,\quad \eta^N=\chi^N=1.
    \end{equation}
    Solving it with respect to $T$ one gets a linear relation for $u,v,w$, which reads 
    \begin{equation}
        u+\chi\cdot  v = \eta \cdot w,
    \end{equation}
    and gives all the linear factors of \eqref{eq:factorLoci}. \newline
    2) $Y_k=-1$. It implies that $z=0$, because of the equation.
    
    The singular locus is invariant under the action of the direct product of Galois groups of $\mathbb{Z}[\exp(2\pi i/N)]$. This implies that resulting locus in $x,y,z$ coordinates is a polynomial of degree $N$ with integer coefficients.
\end{proof}
The theorem above leads to the following statement:
\begin{theorem}
    The polynomials $\Delta_N(x,y,z)$ may be expressed via $T$-discriminants of the $2N-2$ degree polynomials given by
    \begin{equation}
        P_{x,y,z}(T)=x T^{N-1} (1+T)^{N-1} + y(1+T)^{N-1}-T^{N-1}z.
    \end{equation}
    More precisely, the equality holds
    \begin{equation}
        (-1)^N(N-1)^{2(N-1)}(xyz)^{N-2}\Delta_N(x,y,z) = \operatorname{disc}_T(P_{x,y,z}(T))
    \end{equation}
\end{theorem}
\begin{proof}
    Due to the previous proof, one may consider the non-trivial irreducible part $\Delta_N(x,y,z)$ of the singularity locus as a discriminant-type polynomial. Restricting to the diagonal $Y_1=Y_2=\dots =Y_{N-1}$ we obtain an $2(N-1)$-fold covering of $\mathbb{P}^{2}$ given by the polynomial.
    Such a covering is ramified along discriminant divisor of $P_{x,y,z}(T)$ and corresponds to a singular fiber in the original family.
\end{proof}

\begin{figure}
\begin{tikzpicture}[scale=2]
     \draw[very thick](0,-0.15) circle (0.85);
     \draw[very thick,dashed] (2,-1) -- (-2,-1);
     \draw[very thick,dashed] (9/5,-3/2) -- (-1/5,11/6);
     \draw[very thick,dashed] (-9/5,-3/2) -- (1/5,11/6);
     
     \filldraw (0.73, 0.27) circle (1pt) node[anchor=west]{\small $[1:1:0]$};
     \filldraw (-0.73, 0.27) circle (1pt) node[anchor=east]{\small $[0:1:1]$};
     \filldraw (-0, -1) circle (1pt) node[anchor=north]{\small $[1:0:1]$};
\end{tikzpicture}\quad
\begin{tikzpicture}[scale=1.5]
\filldraw (0, 1) circle (2pt) ;
\filldraw (1, 0) circle (2pt) ;
\filldraw (0, -1) circle (2pt) ;
\filldraw (-1, 0) circle (2pt) ;
\draw[dashed, very thick,domain=-2.3:2.3] plot (0,\x);
\draw[dashed, very thick,domain=-2.3:2.3] plot (\x,0);
\draw[very thick,domain=-2:1] plot (\x,\x+1);
\draw[very thick,domain=-2:1] plot (-\x,\x+1);
\draw[very thick,domain=-1:2] plot (\x,\x-1);
\draw[very thick,domain=-1:2] plot (-\x,\x-1);
\end{tikzpicture}
\caption{$N=2$ singularity loci and unfolding in $(u,v)$-chart ($w=1$)}
\end{figure}
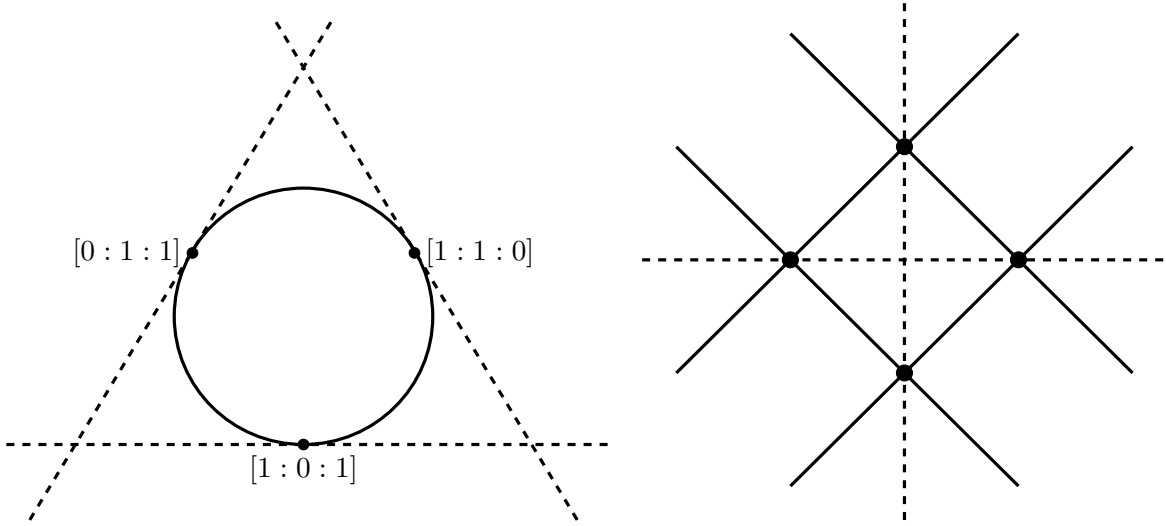

\begin{theorem}\label{th:real}
The projective curve $\Delta_N=0$ is a rational curve that has $(N-1)(N-2)/2$ double points. All singularities belong to $\mathbb{P}^2(\mathbb{R})$. 
Moreover, coordinates of the singularities belong to $\mathbb{Q}(\cos(2\pi/N))$.

\end{theorem}

\begin{proof}
    Singularities of $\Delta_N(x,y,z)$ unfold to singularities of $\Delta_N(u^N,v^N, w^N)$, which are just intersection of it irreducible components (lines). These singularities are given by
    $$
    [u:v:w]=[\chi\xi(1-\alpha\beta) \,:\, -\chi(1-\alpha) \,:\, \xi(1-\beta)], \quad \alpha^N=\beta^N=\chi^N=\xi^N=1.
    $$
    Returning back to initial coordinates $[x:y:z]$, singular points become
    $$
    [(1-\alpha\beta)^N \,:\,  (-1)^N(1-\alpha)^N \,:\, (1-\beta)^N].
    $$
    Now one may see that coordinates are invariant with respect to the Galois conjugation, which sends $\alpha$ to $\alpha^{-1}$ and $\beta$ to $\beta^{-1}$. Indeed, one has that singularities transfer to
    \begin{equation*}
        [(\alpha\beta)^{-N}(\alpha\beta-1)^N \,:\,  (-1)^N\alpha^N(\alpha-1)^N \,:\, \beta^N(\beta-1)^N]
    \end{equation*}
    which is the same as initial coordinate after simplification
\end{proof}

\section{Multiple products and multinomial coefficient}\label{sec:LG2}
The associative structure given by the kernel $K_N(x,y,z)$ allows one to evaluate multiple products of the $N$-Bessel functions via convolution with some higher kernel. This leads to identities
\begin{equation}
    \Phi_N(x_1)\Phi_N(x_2)\dots \Phi_N(x_{m-1})\Phi_N(x_m) = \frac{1}{2 \pi i} \oint K_N(x_1,x_2,\dots x_m \, | z)\,\Phi_N(z) \frac{dz}{z}.
\end{equation}
Obviously, the kernel $K_N(x_1,x_2,\dots x_m \, | z)$ may be obtain by a convolution of the multiplication kernel with itself several times. For the higher kernel one can prove, using the structure constant approach that $K_N(x_1,x_2,\dots x_m \, | z)$ is the generating function for the powers of multinomial coefficients, i.e. it reads
\begin{equation}
    K_N(x_1,x_2,\dots x_m \, | z) = \sum\limits_{r=0}^\infty\sum\limits_{j_1+j_2 + \dots + j_m=r} \binom{r}{j_1,j_2,\dots\, j_m}^N\frac{x_1^{j_1}x_2^{j_2}\dots x_m^{j_m}}{z^{r}}.
\end{equation}
As in the case of the standard product formula, the power series above may be written as the period for a $\mathbb{P}^{m}$-family of hypersurfaces in $\mathbb{G}_{ (N-1)\times (m-1)}$.
\begin{theorem}\label{th:mult}
Consider a $\mathbb{P}^{m}$-family of hypersurfaces in $\mathbb{G}_{ (N-1)\times (m-1)}$ given by the equation
\begin{equation}\label{eq:LGNm}
    W_N^{(m)} (x_1,\dots x_m, z)= \prod\limits_{j=1}^{N-1}\left(1+\sum\limits_{l=1}^{m-1}Y_j^{(l)}\right) \cdot \left(x_1+\sum\limits_{l=1}^{m-1} \prod\limits_{j=1}^{N-1}\frac{x_{l+1}}{Y_j^{(l)}}\right) - z =0.
\end{equation}
Then its period given by the Poincar\'e residue
    \begin{equation}\label{eq:DublPeriod}
        \pi_N(x_1/z,x_2/z, \dots x_m/z) = \frac{1}{(2\pi i)^{N-1}}\oint \frac{1}{W_N^{(m)}(x_1/z,x_2/z,\dots x_m/z,1)} \prod_{j,l}\frac{\operatorname{d} Y_{j}^{(l)}}{Y_{j}^{(l)}}, 
    \end{equation}
    coincides with the formal power series
    $$
   \sum\limits_{r=0}^\infty \sum\limits_{j_1+j_2 + \dots + j_m=r} \binom{r}{j_1,j_2,\dots\, j_m}^N \frac{x_1^{j_1}x_2^{j_2}\dots x_m^{j_m}}{z^{r}},
    $$
    for the small values of $z$.
    \end{theorem}
    \begin{proof}
        The proof is similar to the one for the product formula, but now one should take the function
        \begin{equation*}
            F = \left(\alpha +\sum\limits_{l=1}^{m-1}\beta_i Y^{(l)}_0\right )\prod\limits_{j=1}^{N-1}\left(1+\sum\limits_{l=1}^{m-1} Y^{(l)}_j\right), 
        \end{equation*}
        and restrict it to the submanifold
        \begin{equation*}
            Y^{(k)}_0Y^{(k)}_1\dots Y^{(k)}_{m-1} = 1,\quad k=1\dots m-1.
        \end{equation*}
        Then the rest is similar to Theorem \ref{th:prod}.
    \end{proof}

   \begin{theorem} The singular locus of the Picard-Fuchs system for $K^{(m)}_N$ is given by 
$$
x_0x_1x_2x_3\dots x_m\Delta(x_0,x_1,\dots x_m)=0,
$$
where $\Delta$ is symmetric polynomial given by
\begin{equation}\label{eq:factorLoci}
    \Delta(u_0^N,u_1^N\dots u_m^N)=\prod_{\omega_i}^{\omega_i^N=1}\left(u_0 +\sum\limits_{i=1}^m \omega_iu_i\right),
\end{equation}
where $\omega_i$ runs independently over the $N$-th roots of unity.
\end{theorem}

\section{Duplication formulas: Periods and corresponding differential equations}\label{sec:per}

In this section we study a Picard-Fuchs theory of the generating series given by the Landau-Ginzburg potential $W^{(2)}_N$ restricted to the diagonal. These restricted superpotentials corresponds to the Clausen-type duplication formulas which express $\Phi_N(x)^2$ via $\Phi_N(x)$. One may think of the kernel as a lift to the second symmetric power for $N$-Bessel equation. Such kernel is given by the period of the form
\begin{equation}\label{eq:per}
    \pi_N(t) = \frac{1}{(2\pi i)^{N-1}}\oint\frac{d Y_1}{Y_1}\frac{d Y_2}{Y_2}\dots \frac{d Y_{N-1}}{Y_{N-1}} \frac{1}{1-tV_N} = \sum\limits_{i=0}^\infty c_i t^i,\quad c_i =  [ (V_N)^i]_0.
\end{equation}
where $V_N$ is a Landau-Ginzburg superpotential which reads as
\begin{equation*}
    V_N = \prod\limits_{i=1}^{N-1}(1+Y_i)+\prod\limits_{i=1}^{N-1}(1+Y^{-1}_i) = -W_N(1,1,0)
\end{equation*}
Since we restrict to the diagonal, the base of deformation becomes one dimensional and Gauss-Manin connection transfers to a Picard-Fuchs operator in one variable. Here we provide the first few operators which annihilate $\pi_N(t)$ for different $N$:
\begin{align*}
    N=2\quad & \left(4t - 1\right) \theta_t + 2t,\\
    N=3\quad & \left(t+1\right)(8t-1) \theta_t^{2} + t\left(16t + 7\right) \theta_t + 2t(4t +1),\\
    N=4 \quad & (16t - 1)(4t + 1) \theta_t^{3} + 6 t {\left(32 t + 3\right)} \theta_t^{2} + 2 t {\left(94  t + 5\right)}  \theta_t + 2t(30t + 1), \\
    N=5 \quad & (32t - 1)  (4t - 7)^{2}  (t^{2} - 11t - 1) \theta_t^{4} + 2 t  (4t - 7)  (256t^{3} - 2084t^{2} + 4942t + 143) \theta_t^{3} +\\ &  t  (3072t^{4} - 23024t^{3} + 72568t^{2} - 102261t - 1638) \theta_t^{2} +\\ & t  (2048t^{4} - 12896t^{3} + 30072t^{2} - 66094t - 637) \theta_t + \\ &  2  t  (256t^{4} - 1472t^{3} + 1904t^{2} - 7868t - 49),\\
    N=6 \quad &(t - 1)  (27t + 1)  (64t - 1)  (75t^{3} + 1420t^{2} + 561t + 9) \theta_t^{6} +\\ & +(842400t^{6} + 18022725t^{5} - 1363487t^{4} - 4622791t^{3} - 127551t^{2} - 1977t - 9) \theta_t^{5} +\\ &+ 5  t  (452880t^{5} + 10962507t^{4} + 2491544t^{3} - 1779376t^{2} - 46584t - 168) \theta_t^{4} + \\ &+5  t  (644760t^{5} + 17135271t^{4} + 6994741t^{3} - 1716533t^{2} - 56024t - 96) \theta_t^{3} +\\ & +2  t  (1282500t^{5} + 36696915t^{4} + 19164721t^{3} - 2088858t^{2} - 100236t - 72) \theta_t^{2} + \\ &+2  t  (540900t^{5} + 16436910t^{4} + 9826066t^{3} - 428487t^{2} - 38811t - 9) \theta_t +\\ &+12  t^{2}  (15750t^{4} + 503175t^{3} + 327205t^{2} - 1845t - 1044).
\end{align*}
Here $\theta_t=t\frac{d}{dt}$. In these first six examples, the order of differential operator increases by $1$ each time $N$ increases by $1$ when $N<6$. Then it jumps by $2$ each time we increase $N$ by $1$. The next increasing of slope happens when $N=10$. For $N$ up to $15$ we have
\vspace{0.2cm}
\begin{center}
\begin{tabular}{|c|c|c|c|c|c|c|c|c|c|c|c|c|c|c|c|c|}
    \hline $N$ &2 & 3 & 4 & 5& 6& 7& 8& 9 & 10 & 11 & 12 & 13 & 14 & 15 & 16 & 17\\
    \hline ord$\mathcal{D_N}$ & 1 & 2 &3 & 4& 6& 8&  10& 12& 15 & 18 & 21 & 24 & 28 & 32 & 36 & 40\\
    \hline
\end{tabular}
\end{center}
\vspace{0.2cm}
Our hypothesis is that the slope raises by $1$ on each fifth iteration. However, it has to be checked for larger $N$. All the operators obtained by that procedure are Fuchsian and carry very specific monodromy at $0$. For $2\leq N\leq 5$, the monodromy matrix at $0$ is maximal unipotent with all exponent $=0$, we call such singularities {\em MUM points} (Maximal Unipotent Monodromy point). For $N=6$ the monodromy matrix at $0$ is conjugated to
\begin{equation*}
    M^{(6)}_0 \sim \left(\begin{array}{cccccc}
        1 &1/2! & 1/3! & 1/4! & 1/5! & 0  \\
        0 & 1 &1/2! & 1/3! & 1/4! & 0\\
        0 & 0 & 1 &1/2! & 1/3! & 0 \\
        0& 0& 0 &1 &1/2! & 0\\
        0& 0& 0 & 0 &1 & 0\\
        0& 0& 0 & 0 &0 & 1\\
    \end{array}\right)
\end{equation*}
where the first Jordan block corresponds to the integral analytical solution and corresponding logarithmic solutions, while one-dimensional subspace corresponds to the analytic solution which starts with $t$. For $N=7$ the sizes of both blocks grow by $1$, i.e. we have
\begin{equation*}
    M^{(7)}_0  \sim \left(\begin{array}{cccccccc}
        1 &1/2! & 1/3! & 1/4! & 1/5! &1/6! & 0 &0 \\
        0 & 1 &1/2! & 1/3! & 1/4! &1/5! & 0 &0\\
        0 & 0 & 1 &1/2! & 1/3! & 1/4! &0  &0\\
        0& 0& 0 &1 &1/2! & 1/3! &0 &0\\
        0& 0& 0 & 0 &1 & 1/2! &0 &0\\
        0& 0& 0 & 0 &0 & 1 &0 &0\\
        0& 0& 0 & 0 &0 & 0& 1& 1/2! \\
        0& 0& 0 & 0 &0 & 0& 0& 1 
    \end{array}\right)
\end{equation*}
If our hypothesis is correct for arbitrary $N$, then the number of Jordan blocks equals to $\lceil (N-1)/4\rceil$, where each block $\mathcal{J}_k$ corresponds to the analytic solution which starts with $t^{k-1}$. Even though $0$ is not a MUM point, we still have enough information to do a computation of a 'mirror map' and 'enumerative
invariants' by considering corresponding logarithmic solutions.

The unipotent block of size $N-1$, which corresponds to the analytical solution $\phi_0(t):=\pi_N(t)$, implies that there exists a set of power series  $\phi_i\in t \mathbb{C}[t], \,i=1\dots N-1$, such that 
\begin{equation*}
  y_0(t)=\phi_0(t),\quad 
  y_1(t)=\frac{\phi_0(t)}{2\pi i}\ln(t)+\phi_1(t),\quad y_k(t) = \sum\limits_{j=0}^k\frac{\phi_{k-j}(t)}{(2\pi i )^j}\frac{\ln(t)^j}{j!}, 
\end{equation*}
provide a subset of local solutions corresponding to the unipotent block. Using these solutions, the mirror map coordinate reads as
\begin{equation}
    q = \exp(y_1(t)/y_0(t))\simeq t+O(t^2).
\end{equation}
Usually, there is a suitable choice of the Cauchy data which allows to give an integral mirror map. Inverting this series, one gets a mirror map
\begin{equation}
    t(q) = q+\sum\limits_{j=1}^{\infty} m_j q^j \in \mathbb{Z}[[q]].
\end{equation}
Above we provide the first few coefficients of the expansion for different $N$:
\begin{align*}
    N=3\quad & q-3 q^2+3 q^3+5 q^4-18 q^5+15 q^6+24 q^7-75 q^8+57 q^9+86 q^{10} + O(q^{11})\\
    N=4 \quad & q-4 q^2-6 q^3+56 q^4-45 q^5-360 q^6+894 q^7+960 q^8-6951 q^9+4660 q^{10}+O\left(q^{11}\right) \\
    N=5 \quad & q-5 q^2-40 q^3+115 q^4-645 q^5-12846 q^6-177350 q^7-\\& -2574585 q^8-44198680 q^9-736554815 q^{10}+O\left(q^{11}\right)\\
    N=6 \quad & q-6 q^2-135 q^3-380 q^4-24960 q^5-696366 q^6-\\& -26153302 q^7-901888104 q^8-35369115894 q^9-1381135576280 q^{10}+O\left(q^{11}\right)\\
    N=7 \quad & q-7 q^2-371 q^3-4543 q^4-378637 q^5-20096783 q^6-1568975093 q^7\\&-112310305031 q^8-9251250532328 q^9-758736375700793 q^{10}+O\left(q^{11}\right)
\end{align*}
Knowledge of the mirror map allows us to compute the Yukawa-type coupling as follows
\begin{equation}
    Y(q) = \left(q\frac{d}{dq}\right)^2 \frac{y_2(x(q))}{y_0(x(q))}.
\end{equation}
One of the biggest breakthrough in mirror symmetry was made by P. Candelas, X. de la Ossa, P. Green and L. Parkes \cite{candelas1991exactly}, who conjectured that the Lambert series of the Yukawa coupling is a generating series of the Gromov-Witten invariants for the corresponding mirror Calabi-Yau threefold. In particular, the Yukawa coupling reads
\begin{equation}
    Y(q) = 1 + \sum\limits_{d=1}^\infty n_d d^3 \frac{q^d}{1-q^d},
\end{equation}
where $n_d$ stand for instanton numbers (BPS numbers). In our case this happens only for $N=5$, where one has a Calabi-Yau operator AESZ: 22\footnote{AESZ number is a classification id in the Calabi-Yau differential operator database, \href{https://cydb.mathematik.uni-mainz.de/}{link}}. However, if we drop an assumption that the differential operator is a Picard-Fuchs operator for a family of some Calabi-Yau threefolds, the integrality property of the coefficients for a Lambert series doesn't disappear. However, one needs to modify the usual recipe for the instanton numbers. One finds experimentally that the coefficients are only divisible by $d^2$, whereas the divisibility by $d^3$ is a strong special property of CY threefolds. Moreover, it happens that $n_d$ are only nearly integral and one has to allow denominators involving finitely many primes.

\section{Deformed mirror map: $4$-Bessel example}\label{sec:def}
Restriction to the diagonal leads to a great source of interesting Picard-Fuchs equations and give rise to enumerative invariants as we showed in a previous section. However, a more challenging problem is to study the properties of the whole $\mathbb{P}^2$ family. In the case of the $m=2$ multiplication formulas, one may see the whole family as a special deformation of the diagonal. In this section we provide an example of the $4$-Bessel multiplication kernel. Before we start our exposition, let us perform a change of coordinates in the $z=1$ affine chart of $\mathbb{P}^2$
\begin{equation*}
    x=t,\quad y=rt,\quad z=1.
\end{equation*}
Such change of variables coincides with the restriction to the one-dimensional ray with slope $r$. In this section, we treat $r$ as a deformation parameter, while $t$ will be a base coordinate for the corresponding Gauss-Manin connection. In that case, Landau-Ginzburg potential reads as
\begin{equation}
    V^{(4)} = r(1+X)(1+Y)(1+Z)+(1+X^{-1})(1+Y^{-1})(1+Z^{-1}),
\end{equation}
while the period becomes
\begin{equation*}
    y_0(t;r):=\pi(t;r) = \frac{1}{(2\pi i)^3}\oint \frac{dX}{X}\frac{dY}{Y}\frac{dZ}{Z}\frac{1}{1-tV^{(4)}(X,Y,Z;r)}.
\end{equation*}
As in the case of duplication formula, we recover a Picard-Fuchs equation from the $t$-expansion of the period function.  The corresponding Picard-Fuchs equation is an $r$-family of the fourth order differential operators
\begin{equation}\label{eq:PFr}
    P_4(t,r)\theta_t^4+P_3(t,r)\theta_t^3+P_2(t,r)\theta_t^2+P_1(t,r)\theta_t+P_0(t,r)=0,
\end{equation}
where
\begin{small}
\begin{align*}
    P_4(t,r) :=& (1-4(r+1)t + 2(3 r^{2}-62 r +3)t^2-4 \left(1+r \right) \left(r^{2}+30 r +1\right)t^3+\left(r -1\right)^{4}t^4)\cdot \\ & \cdot(-5+20(r+1)t-2(4 r^{2}+37 r +4)t^2+4 \left(1+r \right) \left(2 r^{2}+11 r +2\right)t^3-\left(r^{2}+4 r +1\right) \left(r -1\right)^{2}t^4),\\
    &\\
    P_3(t,r) :=& 5-10(r+1)t-2(68r^2 - 1311r + 68)t^2+2 \left(1+r \right) \left(253 r^{2}-1756 r +253\right)t^3-\\
    &-2(317 r^{4}+6587 r^{3}+4992 r^{2}+6587 r +317)t^4 +\\&+ 2 \left(1+r \right) \left(133 r^{4}+2880 r^{3}+23254 r^{2}+2880 r +133\right)t^5 +\\ &+ 14(4 r^{6}+279 r^{5}+2032 r^{4}+3050 r^{3}+2032 r^{2}+279 r +4)t^6-\\&  -2 \left(1+r \right) \left(29 r^{4}+254 r^{3}+394 r^{2}+254 r +29\right) \left(r -1\right)^{2}t^7+5 \left(r^{2}+4 r +1\right) \left(r -1\right)^{6}t^8,\\
    &\\
    P_2(t,r) :=& 25(r+1)t-(151 r^{2}-3622 r +151)t^2-\left(1+r \right) \left(583 r^{2}-456 r +583\right)t^3+\\&+(1009 r^{4}+26366 r^{3}+50930 r^{2}+26366 r +1009)t^4-\\& -5 \left(1+r \right) \left(131 r^{4}+1376 r^{3}+11818 r^{2}+1376 r +131\right)t^5+\\& + (r^{6}-4858 r^{5}-37969 r^{4}-50668 r^{3}-37969 r^{2}-4858 r +1)t^6 \\ & 
    \left(1+r \right) \left(111 r^{4}+646 r^{3}+166 r^{2}+646 r +111\right) \left(r -1\right)^{2}t^7-9 \left(r^{2}+4 r +1\right) \left(r -1\right)^{6}t^8,\\
    &\\
    P_1(t,r):=&-5(r+1)t+5(17 r^{2}-454 r +17)t^2-\left(1+r \right) \left(407 r^{2}+2996 r +407\right)t^3+\\&+(815 r^{4}+21482 r^{3}+52566 r^{2}+21482 r +815)t^4-\\
    &-\left(1+r \right) \left(643 r^{4}+3624 r^{3}+33346 r^{2}+3624 r +643\right)t^5\\ & (67 r^{6}-2394 r^{5}-21883 r^{4}-25500 r^{3}-21883 r^{2}-2394 r +67)t^6 +\\& +5 \left(1+r \right) \left(19 r^{4}+82 r^{3}-82 r^{2}+82 r +19\right) \left(r -1\right)^{2}t^7-7 \left(r^{2}+4 r +1\right) \left(r -1\right)^{6}t^8,\\
    &\\
    P_0(t,r) :=&-20\left( r^{2}-26 r +1\right) t^{2}+10 \left(1+r \right) \left(11 r^{2}+146 r +11\right) t^{3}-\\&-10 \left(r^{2}+22 r +1\right) \left(25 r^{2}+74 r +25\right)t^{4}+\\&+4 \left(1+r \right) \left(56 r^{4}+173 r^{3}+1732 r^{2}+173 r +56\right) t^{5}+\\&+\left(-36 r^{6}+372 r^{5}+4644 r^{4}+4440 r^{3}+4644 r^{2}+372 r -36\right) t^{6}-\\&-6 \left(1+r \right) \left(5 r^{4}+18 r^{3}-34 r^{2}+18 r +5\right) \left(r -1\right)^{2} t^{7}+2 \left(r^{2}+4 r +1\right) \left(r -1\right)^{6} t^{8}.
\end{align*}
\end{small}
The first factor in $P_4(t,r)$ defines the set of non-apparent singularities. The monodromy for the Picard-Fuchs equation has a unipotent block of size $3$ which corresponds to the analytic solution. This data allows us to define both a mirror map for the $r$-deformed family, as well as $r$-deformed 'enumerative invariants'.
Two satellite logarithmic solutions may be chosen such as
\begin{align*}
    y_0(t;r) := &\pi(t;r),\\
    y_1(t;r) := &\ln(t) y_0 + \phi_1(t;r),\\  y_1(t;r) :=&\frac{\ln(t)^2}{2}y_0 + \ln(t)\phi_1(t;r)+\phi_2(t;r),
\end{align*}
with
\begin{align*}
\phi_1(t;r):=&2\left(r+1 \right) t +\left(3 r^{2}+32 r +3\right) t^{2}+\frac{r+1}{3}\left(11r^{2}+556 r +11\right) t^{3}+\mathrm{O}\! \left(t^{4}\right),\\
\phi_2(t;r):=&-\frac{1}{5}\left(7 r^{2}+ r +7\right) t^{2}-\frac{r+1}{15}  \left(7r^{2}+346 r +7\right) t^{3}+\\&+\frac{1}{6}\left(3r^{4}-498 r^{3}+1634 r^{2}-498 r +3\right) t^{4}+\mathrm{O}\! \left(t^{5}\right).
\end{align*}
Resolving the mirror map equation, we get a $q$-expansion for the initial coordinate 
\begin{align*}
    t(q;r)=& q-2(r+1)q^2+(5 r^{2}-16 r +5)q^3-14 \left(1+r \right) \left(r^{2}-4 r +1\right)q^4+\\&+(42 r^{4}-108 r^{3}+87 r^{2}-108 r +42)q^5 -\\& -4 \left(1+r \right) \left(33 r^{4}-104 r^{3}+187 r^{2}-104 r +33\right)q^6 +\mathrm{O}\! \left(q^{7}\right).
\end{align*} 
Computation of the Yukawa potential gives the Lambert series
\begin{equation}
    Y(q;r)=1+\sum\limits_{i=1}^{\infty} a_i\frac{q^i}{1-q^i},
\end{equation}
where the instanton-type numbers are
\begin{align*}
    \, & a_{1} = 0, \quad a_{2} = -4 \left(r -1\right)^{2}
,\quad  a_{3} = 18 \left(1+r \right) \left(r -1\right)^{2} \\ &
a_{4} = 
-\frac{2 \left(107 r^{2}+182 r +104\right) \left(r -1\right)^{2}}{3}, \quad 
a_{5} = 
5 \left(1+r \right) \left(55 r^{2}+64 r +55\right) \left(r -1\right)^{2},\\
& a_{6} = 
-\frac{\left(15797 r^{4}+40042 r^{3}+45282 r^{2}+40177 r +15912\right) \left(r -1\right)^{2}}{15}
,\\ & a_{7} = 
\frac{14 \left(1+r \right) \left(4319 r^{4}+8044 r^{3}+6849 r^{2}+8044 r +4319\right) \left(r -1\right)^{2}}{15}.
\end{align*}
While integrality properties are not clear for the $r$-family, every instanton number is divisible by $(r-1)^2$. Such factor indicates that for $r=1$ we have a special behavior which corresponds to the pure duplication formula. We know that in such a case, the Picard-Fuchs operator has to be a third order operator , appeared in  section \ref{sec:per} $N=4$. In the $r$-parametric story this "drop" of order appears as a reducibility of operator \eqref{eq:PFr} - for $r=1$ operator splits into the product of the first order differential operator and the non-trivial third order operator with MUM point at $0$ mentioned above.

\begin{remark}
The multiplication kernels become trivial in the case when $r=0$. However, the corresponding Picard-Fuchs equation remains still nontrivial. The same also holds for the mirror coordinate and other invariants. In particular, one can see that mirror coordinate for $r=0$ reads as
\begin{equation*}
    t(q;0)=q -2 q^{2}+5 q^{3}-14 q^{4}+42 q^{5}-132 q^{6}+429 q^{7}-1430 q^{8}+4862 q^{9}+\mathrm{O}\! \left(q^{10}\right),
\end{equation*}
where one can recognize the Catalan numbers OEIS:\href{https://oeis.org/A000108}{A000108} 
\end{remark}

\section{Buchstaber-Rees $N$-valued groups, $N$-Bessel kernels and $N$-fold singular covers}\label{sec:Buch}
If $X$ is a topological space, let ${\rm Sym}^n (X)$ denote its $n-$fold symmetric product, i.e., ${\rm Sym}^n (X)  = X\times X\times\ldots\times X/\mathfrak{S}_n$ 
where the symmetric group $ \mathfrak{S}_n$  acts on $X^n :=  X\times X\times\ldots\times X$ by permuting the factors.
An $n-$multiset of a given space $X$ is an unordered set of $n$ points of $X$ (multiplicity including).  The symmetric product ${\rm Sym}^n (X)$ is the topological space of all $n-$multisets of $X.$

\begin{definition}
An $n-$valued multiplication on $X$ is a map
\begin{equation*}
    \mu: X\times X \to {\rm Sym}^n (X)
\end{equation*}
given by
\begin{equation}
\mu(x,y)= x\star y =[z_1,z_2,\ldots,z_n], \quad z_k = (x\ast y)_k\end{equation}
\end{definition}
\begin{definition}
A set $X$ with an $n$-valued multiplication is called an $n$-valued group if the following conditions hold:
\begin{itemize}
\item{\bf Associativity:} The following two $n^2-$ sets are equal:
$$[x\star(y\star z)_1,\dots,x\star(y\star z)_n] = [(x\star y)_1\star z,\dots, (x\star y)_n\star z]$$
for all $x; y; z\in X$. This condition is equivalent to the commutative diagram
\begin{center}
\begin{tikzcd}
    & X \times {\rm Sym}^n(X) \arrow[r, "D\otimes 1"] & {\rm Sym}^n(X\times X) \arrow[dr, "\mu^n"] \\
    X\times X \times X\arrow[ur, "1\otimes \mu"] \arrow[dr, "\mu\otimes 1"'] & & & {\rm Sym}^{n^2}(X)  \\
    & {\rm Sym}^n(X) \times X \arrow[r, "1\otimes D"] & {\rm Sym}^n(X\times X) \arrow[ur, "\mu^n"']
\end{tikzcd}
\end{center}
where $1\otimes D([x_1,x_2,\dots x_n],x_0) = [(x_1,x_0),(x_2,x_0),\dots (x_n, x_0)] $ and
$D\otimes 1$ defined in the similar way.
\item{\bf  Unit:} There exists an element $e\in X$ such that 
$$
e\star x = x\star e =[x,x,\ldots, x]
$$
for all $x\in X.$

\item{\bf Inverse element:} There is a map ${\rm inv} : X\to X$ such that the unit $e$ belongs to sets ${\rm inv}(x)\star x =x\star{\rm inv}(x)$
for all $x\in X.$
\end{itemize}
\end{definition}

\subsection{Additive $N$-valued group structure on $\mathbb{C}$ \cite{buchstaber1997multivalued}}
Define a multi-valued binary operation
\begin{equation}
\mu : \mathbb{C}\times \mathbb{C} \to {\rm Sym}^N (\mathbb{C})
\end{equation}
by the formula 
\begin{equation}
x \star y =\left[\left(x^{1/N}+ \chi^r y^{1/N}\right)^{n},\quad 1\leq r\leq N\right]
\end{equation}
where $\chi$ is a primitive $N$-th root of unity. The introduced multiplication endows $\mathbb{C}$ with the structure of an $N-$valued group with the unit
$e = 0.$ The inverse element is given by the map ${\rm inv}(x) = (-1)^Nx$. One can define a divisor which is supported by the $N$-valued multiplication locus as
\begin{equation*}
    \Delta(x,y,z) = (z-z_1)(z-z_2)\dots (z-z_N),\quad \mu(x,y)=[z_1,z_2,\dots z_n].
\end{equation*}
The form of the group law and Theorem \ref{th:sing1} leads to the following result
\begin{theorem}
    The kernel for the $N$-Bessel product formula lifts the Buchstaber-Rees $N$-valued group laws.
\end{theorem}
The importance of this observation underlines the following idea: multiplication kernel lifts the "multiplication" loci from the corresponding spectral curve. By lifting we mean a Green-like function, which replaces the $\delta$-function by the measure with singularity locus concentrated at the support of corresponding $\delta$. Indeed, the one consider the corresponding connection 
$$
\delta = x\frac{d}{dx} - A = x\frac{d}{dx}-\begin{pmatrix}
    0 & 1 & 0 & \dots & 0 & 0 \\
    0 & 0 & 1 & \dots & 0 & 0\\
    \dots & \dots & \dots  & \dots & \dots& \dots \\
    0 & 0 & 0 & \dots & 1 & 0 \\
    0 & 0 & 0 & \dots & 0 & 1 \\
    x & 0 & 0 & \dots & 0 & 0
\end{pmatrix}
$$
The spectral curve then reads 
\begin{equation*}
    \Gamma(x,\lambda)=(-1)^N\operatorname{det}(A/x-\lambda)=\lambda^N-1/x=0.
\end{equation*}
One can see that it corresponds to an $N$-fold cover of $\mathbb{C}$ and may be seen as $\operatorname{Sym}^N(\mathbb{C})$ for $x\in \mathbb{C}^{*}$, which naturally carries the multivalued group structure described above. 

\section{Generalized Frobenius method. Bessel Kernels and symmetric polynomials}\label{sec:FrobNum}

In the previous section we used the special construction of multiplication kernels, which prescribes the expansion as a multi-Taylor series in $x,y,1/z$. In this section we introduce another approach for finding a different type
of expansion of the kernel function. While in the previous part we relate kernels to the periods for some explicit families of algebraic hypersurfaces, here we concentrate on the integrality properties of the introduced expansion. 

\subsection{Generalized Frobenius method}Here we use some special property of the Bessel kernel as well as a general property of the multiplication kernels. Firstly, we use the fact that $N$-Bessel kernels are symmetric functions in $x,y,z$ variables, since the $N$-Bessel operator $\frac{1}{x}\theta^N$ is self-adjoint. Secondly, for any kernel function one may introduce the following Cauchy condition \cite{kontsevich2021multiplication}:
\begin{equation}\label{eq:bound}
  K_{N}(0,y,z)=\frac{1}{z-y}.
\end{equation}
Moreover, the kernel solves a system of partial differential equations
\begin{equation}\label{eq:FrobSys}
    \left[\frac{d}{dx}\theta_x^{N-1}-\frac{d}{dz}\theta_z^{N-1}\right]K=0,\quad \left[\frac{d}{dy}\theta_y^{N-1}-\frac{d}{dz}\theta_z^{N-1}\right]K=0.
\end{equation}
Now we are looking for a solution as a Taylor series in $x$ satisfying "initial" condition \eqref{eq:bound}. Then the kernel reads
\begin{equation}\label{eq:Frob}
   K_{N}(x,y,z) = \sum\limits_{j=0}^{\infty}a_j(y,z)x^j,\quad a_0(y,z)= \frac{1}{z-y},
\end{equation}
and the following two differential recursions hold
\begin{equation}
    a_j(y,z)=\frac{1}{j^N}\frac{d}{dz}\theta_z^{N-1}[ a_{j-1}(y,z)],\quad a_j(y,z)=\frac{1}{j^N}\frac{d}{dy}\theta_y^{N-1}[ a_{j-1}(y,z)].
\end{equation}
The existence of two recursions of the same form for this example proves that the functions  $a_j(y,z)$ are symmetric (up to the sign) under interchange of $y$ and $z$. This formulae allows to rewrite a kernel function for a generic choice of $x$ and $y$. Assuming that an initial point for the recursion is given by $1/(z-y)$, it is easy to see that $K^{(N)}(x,y,z) \in \mathbb{Q}\left[y,z,(z-y)^{-1}\right][[x]]$. We call series \eqref{eq:Frob} a generalized Frobenius solution of \eqref{eq:bound}-\eqref{eq:FrobSys}.

\subsection{Palindromic polynomials and Bessel kernels} The main observation of this section is, that this type of series representation for the kernel function enjoys some integrality properties. Write the expansion as
\begin{equation}\label{eq:FrExp}
   K_{N}(x,y,z) = \sum\limits_{m=0}^{\infty}P^{(N)}_m(y,z) \frac{x^m}{(z-y)^{Nm+1}},\quad P_m(y,z) \in \mathbb{Z}[y,z].
\end{equation}
The polynomials $P^{(N)}_m(y,z)$ are symmetric, i.e.
\begin{equation}
    P^{(N)}_m(y,z) = P^{(N)}_m(z,y),
\end{equation}
so their coefficients define palindromic polynomial. For example, for $N=2$ the kernel reads,
  \begin{multline}
K_2(x,y;z)\simeq \frac{1}{z -y}+\frac{\left(y +z \right)}{\left(z -y \right)^{3}}x+\frac{\left(y^{2}+4 y z+ z^{2}\right)}{\left(z -y \right)^{5}}x^2+\\+\frac{\left(y^{3}+9 y^{2} z +9 y z^{2}+z^{3}\right)}{\left(z -y \right)^{7}}x^3+\frac{\left(y^{4}+16 y^{3} z +36 y^{2} z^{2}+16 y z^{3}+z^{4}\right)}{\left(z -y \right)^{9}} x^{4}+O(x^5).
\end{multline}
The palindromic polynomials $P^{(2)}_i(y,z)$ are homogeneous of degree $i$ with integer coefficients,
$$
P_i^{(2)}(y,z) = \sum\limits_{k+n=i}T^{(i)}_{k,n}y^k z^n,\quad T^{(i)}_{k,n}=T^{(i)}_{n,k}\in\mathbb{Z},\quad T^{(i)}_{i,0}=T^{(i)}_{0,i}=1.
$$ 
Writing down the coefficients $T^{(i)}_{k,n}$ in a triangle form we get
$$
\arraycolsep=0.7pt\def\arraystretch{0.8}
\begin{array}{cccccccccccccc}
&&&&&&&T^{(0)}_{0,0}&&&&&&\\
&&&&&&T^{(1)}_{1,0}&&T^{(1)}_{1,0}&&&&&\\
&&&&&T^{(2)}_{2,0}&&T^{(2)}_{1,1}&&T^{(2)}_{0,2}&&&&\\
&&&&&&&\dots&&&&&&
\end{array}=
\arraycolsep=0.7pt\def\arraystretch{0.8}
\begin{array}{cccccccccccccccc}
&&&&&&&&1&&&&&&&\\
&&&&&&&1&&1&&&&&&\\
&&&&&&1&&4&&1&&&&&\\
&&&&&1&&9&&9&&1&&&&\\
&&&&1&&16&&36&&16&&\,1&&&\\
&&&&&&&&\dots&&&&&&&
\end{array}
$$
where one may recognize the square of the Pascal triangle entries. For $N=3$ and $N=4$ these triangles exist at OEIS: \href{https://oeis.org/A181544}{A181544} for $N=3$,
\href{https://oeis.org/A262014}{A262014} for $N=4$. 

\subsection{Combinatorial properties}The way we obtain such polynomials is very similar to a construction introduced in \cite{agapito2014symmetric} for the Euler-like and Naryana-like polynomials. In our case the polynomials seem to enjoy very similar combinatorial properties. Firstly, one easily gets
\begin{equation*}
    P^{(N)}_m(y,z) = y^{m(N-1)}P^{(N)}_m(z/y,1).
\end{equation*}
The polynomials $P^{(N)}_m(t,1)$ are palindromic polynomials with integer coefficients. Now let us write them as
\begin{equation}
    P^{(N)}_m(t,1) = \sum\limits_{k=0}^{(Nm+1)/2} \gamma_{k}^{(N,m)}t^k(1+t)^{Nm+1-2k}.
\end{equation}
We call the tuple consisting of $(\gamma_{0}^{(N,m)},\gamma_{1}^{(N,m)},\dots )$ a $\gamma$-vector of polynomial $P^{(N)}_m$. Such $\gamma$-vectors are invariants of symmetric reciprocal polynomials. 
\vspace{0.3cm}
\newline
{\bf Conjecture 1.} {\em Entries of the $\gamma$-vectors for the polynomials $P^{(N)}_m$ are positive.}
\vspace{0.3cm}
\newline
In order to give some evidence of this conjecture we record first six nontrivial $\gamma$-vectors for $N=2,3$ and $4$. This brings us to a stronger conjecture
\vspace{0.3cm}
\newline
{\bf Conjecture 2.} {\em All roots of $P^{(N)}_m$ are real. This implies Conj. 1}
\vspace{0.3cm}

\renewcommand{\arraystretch}{1.2}
\begin{table}[]
    
    \begin{tabular}{|l|l|} \hline $N=2$& $y \,+\, z$ \\
         & $y^2 \,+\, 4yz \,+\, z^2$\\
         & $y^3 \,+\, 9y^2z \,+\, 9yz^2 \,+\, z^3$ \\
         & $y^4 \,+\, 16y^3z \,+\, 36y^2z^2 \,+\, 16yz^3 \,+\, z^4$ \\
         & $y^5 \,+\, 25y^4z \,+\, 100y^3z^2 \,+\, 100y^2z^3 \,+\, 25yz^4 \,+\, z^5$ \\
         & $y^6 \,+\, 36y^5z \,+\, 225y^4z^2 \,+\, 400y^3z^3 \,+\, 225y^2z^4 \,+\, 36yz^5 \,+\, z^6$ \\
         \hline
        $N=3$ & $y^2 \,+\, 4yz \,+\, z^2$ \\
& $y^4 \,+\, 20y^3z \,+\, 48y^2z^2 \,+\, 20yz^3 \,+\, z^4$ \\
& $y^6 \,+\, 54y^5z \,+\, 405y^4z^2 \,+\, 760y^3z^3 \,+\, 405y^2z^4 \,+\, 54yz^5 \,+\, z^6$\\
& $y^8 \,+\, 112y^7z \,+\, 1828y^6z^2 \,+\, 8464y^5z^3 \,+\, 13840y^4z^4 \,+\, 8464y^3z^5 \,+\, $\\ &$1828y^2z^6 \,+\, 112yz^7 \,+\, z^8$ \\
& $y^{10} \,+\, 200y^9z \,+\, 5925y^8z^2 \,+\, 52800y^7z^3 \,+\, 182700y^6z^4 \,+\, 273504y^5z^5 \,+\, $ \\ & $182700y^4z^6 \,+\, 52800y^3z^7 \,+\, 5925y^2z^8 \,+\, 200yz^9 \,+\, z^{10}$\\\hline
    \end{tabular}
    \vspace{0.15cm}
    \caption{First few polynomials $P^N_m$}
    \label{tab:my_label}
\end{table}
\renewcommand{\arraystretch}{1.2}
\begin{table}[]
    
    \begin{tabular}{|l|l|l|l|}
          \hline $m$ & $N=2$ & $N=3$  & $N=4$ \\
          \hline
     $2$ &  $(1,\,2)$ & $(1,\,2)$ & $(1,\,8)$\\
     $3$ & $(1,\,6)$
  & $(1,\,16,\,10)$  & $(1,\,66,\,324,\,104)$\\
     $4$ & $(1,\,12,\,6)$  & $(1,\,48,\,198,\,56 )$ & $(1,\,234,\,5076,\,19304,\,11136)$\\ 
     $5$ & $(1,\,20,\,30)$ & $(1,\,104,\,1176,\,2144,\,346)$ & \dots  \\
     $6$ & $(1,\,30,\,90,\,20)$ & $(1,\,190,\,4360,\,21200,\,21650,\,2252)$ & \dots \\
         \hline 
    \end{tabular}
    \vspace{0.15cm}
    \caption{$\gamma$-vectors of $P^{N}_m$}
    \label{tab:my_label}
\end{table}



\subsection{Limit to the Dwork family}
Another interesting feature of these triangles is that the sum of elements of the $m$-th row of the $N$-th triangle gives generalized Catalan numbers $(mN)!/(m!)^N$, i.e.
$$
\mathcal{DW}(N;t)=\sum\limits_{m=0}^{\infty}\left[\sum\limits_{k=0}^{m}{}^{(N)}T^{(m)}_{k,m-k} \right] t^m = \sum\limits_{m=0}\frac{(Nm)!}{(m!)^N}t^m.
$$
One can observe that $\pi(N;t)$ is a period for the Dwork hypersurfaces family in $\mathbb{P}^{N-1}$ given by 
\begin{equation}
    X_1^{N}+X_2^{N}+\dots X_{N}^{N}+t X_1X_2\dots X_N = 0.
\end{equation}
It follows that $\pi(N;t)$ arises as a special limit of the kernel $K_N(x,y;z)$. Let us briefly explain this. 

Indeed, the kernel is homogeneous function, we may consider the following series
$$
(z-y)K^{(N)}(x,y,z) = F(x/z,y/z)=F(X,Y)=\sum\limits_{m=0}^{\infty}P(N)_m(Y,1)\frac{X^m}{(1-Y)^{mN}},
$$
where
$$
P(N)_m(Y,1)=\sum\limits_{k=0}^{n(N-1)}{}^{(N)}T^{(m)}_{k,m(N-1)+1-k}Y^k.
$$
Making the change of variables
\begin{equation*}
    X= (1-Y)^Nt
\end{equation*}
and taking the formal limit as $Y$ approaches to $1$,
$$
\lim_{Y\rightarrow 1}F((1-Y)^Nt, Y)=\mathcal{DW}(N;t),
$$
we get the Dwork period.

\section{Picard-Fuchs equations in three variables}\label{sec:PF}

The main feature of the expansion \eqref{eq:FrExp} is that instead of a multi-Taylor expansion in variables
$x,y,z^{-1}$, it provides an expansion in variable $x$ with rational functions of $y,z$ as coefficients.
This allows us to consider the Picard-Fuchs operator as an $(y,z)$-parametrized family of differential
operators in the $x$ variable. Below we present these operators for $N=2,3$ and $4$. 

\subsection{\boldmath ${ N=2}$ Sonine--Gegenbauer kernel}
The Picard-Fuchs equation takes form
\begin{equation}
\left(x^{2}-2 x y -2 x z +y^{2}-2 y z +z^{2}\right)\frac{d}{d x}K +\left(x -z -y \right) K=0.
\end{equation}
Solving this equation with the initial condition
$$
K \! \left(0,y,z\right) = \frac{1}{z -y}
$$
one gets a immediately the Sonine--Gegenbauer kernel function
\begin{equation}
    K_2(x,y,z) = \frac{1}{\sqrt{x^2+y^2+z^2-2(xy+xz+yz)}}.
\end{equation}
This form of the kernel was expected. More nontrivial examples appear for higher $N$.

\subsection{{\boldmath  ${ N=3}$} lift of Legendre family} The Picard-Fuchs family reads
\begin{equation*}
\begin{gathered}
x \left(2 x -y +z \right) \left(x^{3}+3 x^{2} y -3 x^{2} z +3 y^{2} x +21 x y z +3 z^{2} x +y^{3}-3 y^{2} z +3 z^{2} y -z^{3}\right)\frac{d^2}{dx^2}K+\\
+\big[6 x^{4}+8 x^{3} y -8 x^{3} z -3 x^{2} y^{2}+60 x^{2} y z -3 x^{2} z^{2}-6 y^{3} x -36 x y^{2} z +36 x y z^{2}+6 z^{3} x -y^{4}+\\+4 y^{3} z -6 y^{2} z^{2}+4 y z^{3}-z^{4}\big]\frac{d}{dx}K+\left(x +y -z \right) \left(2 x^{2}-2 x y +2 x z -y^{2}-4 y z -z^{2}\right)K=0.
\end{gathered}
\end{equation*}
A solution of this equation with the usual type Cauchy data
$$
K(0, y, z) = \frac{1}{z-y}
$$
leads to an explicit formula for the kernel function:
\begin{equation}
K_3(x,y,z)=-\frac{{}_{2}^{}{{F_{1}^{}}} \left(\frac{1}{3},\frac{2}{3};1;-\frac{27 x y z}{\left(x +y -z \right)^{3}}\right)}{x +y -z}.
\end{equation}
Here again we obtain an expected result. Indeed, the corresponding Landau-Ginzburg variety reads
\begin{equation*}
    xXY(1+X)(1+Y)+y(1+X)(1+Y)-zXY=0,
\end{equation*}
which defines a special family of complex curves of degree $4$, that is an $\mathbb{A}^3$-family of elliptic curves. In fact, the obtained kernel is a lift of the Gauss period for the Legendre family of elliptic curves $Y^2=X(X-1)(X-t)$, where 
$$t=-\frac{27 x y z}{\left(x +y -z \right)^{3}}.$$
The singular point $t=0$ lifts to the coordinate triangle $xyz=0$, while $t=1$ coincides with the singular locus for the $3$-Bessel kernel. Corresponding mirror map reads
\begin{multline}
    x(q) = q+3\frac{\left(y^2+3 y z+z^2\right)}{(y-z)^3}q^2+3\frac{\left(4 y^4+14 y^3 z+21 y^2 z^2+14 y z^3+4 z^4\right)}{(y-z)^6}q^3+\\+\frac{\left(55 y^6+210 y^5 z+366 y^4 z^2+417 y^3 z^3+366 y^2 z^4+210 y z^5+55 z^6\right)}{(y-z)^9}q^4+\\+3\frac{\left(91 y^8+372 y^7 z+700 y^6 z^2+884 y^5 z^3+924 y^4 z^4+884 y^3 z^5+700 y^2 z^6+372 y z^7+91 z^8\right)}{(y-z)^{12}}q^5+\dots 
\end{multline}
The family itself gives a deformation of Beauville family of type IV \cite{beauville1982familles,zagier2016arithmetic} with singular fibres of type $I_6,I_3,I_2,I_1$. Such a deformation reads
   \begin{equation*}
       (\tilde{X}+\tilde{Y})(\tilde{Y}+\tilde{Z})(x\, \tilde{Z}+y\,\tilde{X})+ z \tilde{X}\tilde{Y}\tilde{Z}=0.
   \end{equation*}
   From that point of view, Beauville case is an isomondromic coalescence of $2$ singularities in the $\mathbb{P}^2$-family described above. One may think that such deformation decomposes one of the
   singular fibres into two singular fibres, such that the product of new monodromies equals to the monodromy of the $I_6$ fiber in the original Beauville family. In particular, it means that the corresponding to the Beauville case monodromy subgroup gives a representation of the congruence subgroup $\Gamma_0(6)$.
\subsection{{\boldmath  ${ N=4}$} } The Picard-Fuchs operators family reads as
\begin{equation*}
    \left(P_4(x,y,z)\frac{d^4}{dx^4}+\sum\limits_{i=0}^3P_i(x,y,z)\frac{d^i}{dx^i}\right)K(x,y,z)=0
\end{equation*}
where 
\begin{small}
\begin{align*}
P_4&:=x^{2} (5 x^{4}-20 x^{3} y -20 x^{3} z +8 x^{2} y^{2}+74 x^{2} y z +8 x^{2} z^{2}+8 y^{3} x +52 x y^{2} z +52 x y z^{2}+\\& +8 z^{3} x -y^{4}-2 y^{3} z+6 y^{2} z^{2}-2 y z^{3}-z^{4}) (x^{4}-4 x^{3} y -4 x^{3} z +6 x^{2} y^{2}-124 x^{2} y z \\& +6 x^{2} z^{2}-4 y^{3} x -124 x y^{2} z -124 x y z^{2}-4 z^{3} x +y^{4}-4 y^{3} z +6 y^{2} z^{2}-4 y z^{3}+z^{4}),\\
&\\
P_3&:=x (55 x^{8}-410 x^{7} y -410 x^{7} z +1044 x^{6} y^{2}-1238 x^{6} y z +1044 x^{6} z^{2}-1134 x^{5} y^{3}+\\&+11634 x^{5} y^{2} z +11634 x^{5} y z^{2}-1134 x^{5} z^{3}+366 x^{4} y^{4}+4326 x^{4} y^{3} z -54984 x^{4} y^{2} z^{2}+\\&+4326 x^{4} y z^{3}+366 x^{4} z^{4}+266 x^{3} y^{5}-12934 x^{3} y^{4} z -111172 x^{3} y^{3} z^{2}-111172 x^{3} y^{2} z^{3}-\\&-12934 x^{3} y z^{4}+266 x^{3} z^{5}-244 x^{2} y^{6}-6554 x^{2} y^{5} z -47652 x^{2} y^{4} z^{2}-90780 x^{2} y^{3} z^{3}-\\&-47652 x^{2} y^{2} z^{4}-6554 x^{2} y z^{5}-244 x^{2} z^{6}+62 x y^{7}+1070 x y^{6} z +2178 x y^{5} z^{2}-3310 x y^{4} z^{3}-\\&-3310 x y^{3} z^{4}+2178 x y^{2} z^{5}+1070 x y z^{6}+62 x z^{7}-5 y^{8}+10 y^{7} z +40 y^{6} z^{2}-170 y^{5} z^{3}+\\&+250 y^{4} z^{4}-170 y^{3} z^{5}+40 y^{2} z^{6}+10 y z^{7}-5 z^{8}),
\\
&\\
P_2&:=155 x^{8}-1085 x^{7} y -1085 x^{7} z +2285 x^{6} y^{2}+2460 x^{6} y z +2285 x^{6} z^{2}-1647 x^{5} y^{3}+18437 x^{5} y^{2} z +\\&+18437 x^{5} y z^{2}-1647 x^{5} z^{3}-295 x^{4} y^{4}-8928 x^{4} y^{3} z -121474 x^{4} y^{2} z^{2}-8928 x^{4} y z^{3}-\\&-295 x^{4} z^{4}+941 x^{3} y^{5}-18779 x^{3} y^{4} z -160962 x^{3} y^{3} z^{2}-160962 x^{3} y^{2} z^{3}-18779 x^{3} y z^{4}+941 x^{3} z^{5}-\\&-413 x^{2} y^{6}-7572 x^{2} y^{5} z -57531 x^{2} y^{4} z^{2}-128168 x^{2} y^{3} z^{3}-57531 x^{2} y^{2} z^{4}-7572 x^{2} y z^{5}-413 x^{2} z^{6}+\\&+63 x y^{7}+1635 x y^{6} z +4077 x y^{5} z^{2}-5775 x y^{4} z^{3}-5775 x y^{3} z^{4}+4077 x y^{2} z^{5}+\\&+1635 x y z^{6}+63 x z^{7}-4 y^{8}+8 y^{7} z +32 y^{6} z^{2}-136 y^{5} z^{3}+200 y^{4} z^{4}-136 y^{3} z^{5}+32 y^{2} z^{6}+8 y z^{7}-4 z^{8},\\
& \\
    P_1&=110 x^{7}-740 x^{6} y -740 x^{6} z +1186 x^{5} y^{2}+3968 x^{5} y z +1186 x^{5} z^{2}-260 x^{4} y^{3}+3940 x^{4} y^{2} z +\\&+3940 x^{4} y z^{2}-260 x^{4} z^{3}-726 x^{3} y^{4}-8352 x^{3} y^{3} z -37164 x^{3} y^{2} z^{2}-8352 x^{3} y z^{3}-726 x^{3} z^{4}+540 x^{2} y^{5}-\\&-4596 x^{2} y^{4} z -40224 x^{2} y^{3} z^{2}-40224 x^{2} y^{2} z^{3}-4596 x^{2} y z^{4}+540 x^{2} z^{5}-122 x y^{6}-528 x y^{5} z -\\&-7038 x y^{4} z^{2}-27824 x y^{3} z^{3}-7038 x y^{2} z^{4}-528 x y z^{5}-122 x z^{6}+12 y^{7}+420 y^{6} z+\\& +1188 y^{5} z^{2}-1620 y^{4} z^{3}-1620 y^{3} z^{4}+1188 y^{2} z^{5}+420 y z^{6}+12 z^{7},\\
& \\
    P_0&=10 x^{6}-70 x^{5} y -70 x^{5} z +68 x^{4} y^{2}+364 x^{4} y z +68 x^{4} z^{2}+56 x^{3} y^{3}+164 x^{3} y^{2} z+\\& +164 x^{3} y z^{2}+56 x^{3} z^{3}-90 x^{2} y^{4}-300 x^{2} y^{3} z +660 x^{2} y^{2} z^{2}-300 x^{2} y z^{3}-90 x^{2} z^{4}+\\&+30 x y^{5}-54 x y^{4} z -696 x y^{3} z^{2}-696 x y^{2} z^{3}-54 x y z^{4}+30 x z^{5}-\\&-4 y^{6}+24 y^{5} z +108 y^{4} z^{2}-256 y^{3} z^{3}+108 y^{2} z^{4}+24 y z^{5}-4 z^{6}
\end{align*}
\end{small}
\vspace{0.2cm}
The Riemann scheme for this Picard-Fuchs equation is, for generic $y,z$:
$$
\left\{\begin{array}{cccc}
    0 & \infty & \alpha & \beta \\
    \hline 
    0 & 1 & 0& 0 \\
    0 & 1& 1 & 1 \\
    0 & 1& 2 & 2 \\
    1 & 1& 3 & 1/2 
\end{array}\right\},\quad A(\alpha, y,z)=0, \quad \Delta_4(\beta, y,z)=0
$$
where
\begin{equation*}
\begin{gathered}
\Delta_4(x,y,z)=x^{4}-4 x^{3} y -4 x^{3} z +6 x^{2} y^{2}-124 x^{2} y z +6 x^{2} z^{2}-4 y^{3} x -124 x y^{2} z -124 x y z^{2}-\\-4 z^{3} x +y^{4}-4 y^{3} z +6 y^{2} z^{2}-4 y z^{3}+z^{4},
\end{gathered}
\end{equation*}
\begin{equation*}
\begin{gathered}
A(x,y,z) = 5 x^{4}-20 x^{3} y -20 x^{3} z +8 x^{2} y^{2}+74 x^{2} y z +8 x^{2} z^{2}+8 y^{3} x +52 x y^{2} z +52 x y z^{2}+\\+8 z^{3} x -y^{4}-2 y^{3} z +6 y^{2} z^{2}-2 y z^{3}-z^{4}.
\end{gathered}
\end{equation*}
The zeros of polynomial $A(x,y,z)$ correspond to the apparent singularities of the Picard-Fuchs equation. These apparent singularities appear due to the projection of the three-dimensional family to the one-dimensional $x$-direction. The zeros of $\Delta_4(x,y,z)$ together with $0$ and $\infty$ coincide with the singularity set of the corresponding algebraic family \eqref{eq:LGFam1}, whose projective model in $\mathbb{P}^3$ reads
$$
    xXYZ(1+X)(1+Y)(1+Z)+y(1+X)(1+Y)(1+Z)-zXYZ=0.
$$

The monodromy matrices at $\infty$ and $0$ are conjugated to 
$$
M_\infty\sim \left(\begin{array}{cccc}
    1 & 1 & 0 & 0 \\
    0 & 1 & 1 & 0 \\
    0 & 0 & 1 & 1 \\
    0 & 0 & 0 & 1 
\end{array}\right),\quad M_0\sim \left(\begin{array}{cccc}
    1 & 1 & 0 & 0 \\
    0 & 1 & 1 & 0 \\
    0 & 0 & 1 & 0 \\
    0 & 0 & 0 & 1
\end{array}\right).
$$
Once again, the obtained operator carries the property of the kernel expected before. Indeed the singularity locus $\Delta_4$ coincides with the Buchstaber-Rees polynomial. The appearance of the apparent singularities is a result of projection of three-dimensional local system to a one-dimensional direction.

\section*{Discussion}
As we mentioned in the introduction the study of  kernel type functions for classical Bessel functions is motivated by problems in theoretical physics and special Feynman integrals. The obtained Landau-Ginzburg models and corresponding generating functions were also studied in the framework of Calabi-Yau differential equations and mirror symmetry. Another important connection which we wish to underline here is a relation to the irregular Hodge theory initially established by Deligne \cite{deligne2007singularites}. 

Consider a matrix of so-called Bessel moments, which reads as
\begin{equation*}
    \left(\mathcal{P}_{k}\right)_{ij}=\int\limits_{0}^\infty I_0(z)^i K_0(z)^{k-i} z^{2j-1}dz,\quad 1\leq i,j\leq \lfloor (k-1)/2 \rfloor,
\end{equation*}
together with the {\em Bernoulli matrix} 
\begin{equation*}
    \left(\mathcal{B}_k\right)_{ij}=(-1)^{k-i}\frac{(k-i)!(k-j)!}{k!}\frac{B_{k-i-j+1}}{(k-i-j+1)!},
\end{equation*}
where $B_n$ stands for corresponding Bernoulli numbers. It was conjectured and checked numerically by Broadhurst and Roberts \cite{broadhurst2018quadratic}  for odd moments $k\in 2\mathbb{Z}+1$,, that there exists a matrix $\mathcal{D}_k$ which may be obtained reductively in $k$, such that the following quadratic relation holds
\begin{equation}
    \mathcal{P}_{k}\,\cdot \,\mathcal{D}_k\,\cdot\,\mathcal{P}_{k}=(-2\pi i)^k \mathcal{B}_k.
\end{equation}
This relation was proved by Fres\'an, Sabbah and Yu and enlarged for even moments $k\in 2\mathbb{Z}$ \cite{fresan2023quadratic}. Using the  irregular Hodge theory, they interpret the Bessel moment matrix $\mathcal{P}_k$ as a periods of the Kloosterman motive $M_k$ (introduced in \cite{fresan2022hodge}) with $\mathcal{B}_k$ as a corresponding intersection matrix.

While the formulas we consider correspond to a very specific case of product for $I_0$ Bessel function, the results we obtain should be naturally extended to the case of more complicated products which evolve both $I$ and $K$ type Bessel functions. For example, we recognize looking at Vanhove formula \cite{vanhove2014physics}
\begin{equation*}
    \int\limits_{0}^\infty I_0(z) K_0(z)^{l+1}\,z\,dz = \frac{1}{2^l}\int\limits_{X_i\geq 0} \frac{1}{1-W(X)}\prod\limits_{k=1}^l \frac{dX_k}{X_k},\,\,\, W=\left(1+\sum\limits_{i=1}^{l+1} X_i\right)\left(1+\sum\limits_{i=1}^{l+1} X_i^{-1}\right),
\end{equation*}
the same differential form as in Theorem \ref{th:mult}, by putting $N=2$ and $m=l+1$ in formula \eqref{eq:LGNm}. The difference is that now the integration goes along an other relative cycle. This cycle should correspond to the Stokes ray on the irregular side. Indeed both $I_0$ and $K_0$ are expressed as integrals of the same exponential $1$-form but taken over different relative cycles. An illustrative example is a product formula for $K_0$ functions which reads as \cite{durand1975nicholson,vilenkin1978special}
$$
  K_{0}(x)K_{0}(y)= \int_0^{\infty}\frac{K_0(z)zdz}{\sqrt{x^4+y^4+z^4-2(x^2y^2+x^2z^2+y^2z^2)}},
$$
where once again the same function appears as kernel. 

This allows us to think that the obtained period functions play a universal role and may be used as the periods of the corresponding Kloosterman motives related to the higher Bessel functions.

Recently the authors of \cite{fresan2022hodge} studied the same questions for the higher analogs of Bessel functions in the form of $N$-Kloosterman connection. This generalization naturally fits into our results and should provide the same form of the superpotential.

An open question which we leave for further studies is the following - should we consider the relation between periods of a Landau-Ginzburg model $W$ and Kloosterman motives as an instance of mirror symmetry? We see, similarly to the case of Givental mirror construction \cite{givental1995homological} which connects Picard-Fuchs theory (logarithmic connections/periods) with quantum cohomology (irregular connection/oscillatory integrals),  the transition from irregular connection -- $N$-Bessel equation, to the Fuchsian systems -- Picard-Fuchs for the Bessel kernels.

\bibliographystyle{abbrv}
\bibliography{Bessel.bib}

\newpage

\section*{Appendix: Plotting}

Using a standard computational drawing routines, one may easily miss singularities of $\Delta_N(x,y,z)$, which are just isolated point on the real plane. However, the standard way to see such points is to deform equation $\Delta_N(x,y,z)=0$ to $\Delta_N(x,y,z)=\epsilon$, where $\epsilon$ is a real number. Due to theorem \ref{th:real}, the singularities are real points, since that restricting to the affine chart $z=1$ and such $\epsilon$-deformation will deform these isolated points into small ovals. The connected component looks very similar to the $x,y$-axis, however for odd $N$ it intersects both at points $(1,0)$ and $(0,1)$, and for even $N$ it is tangent to both axis at these points. Above we demonstrate this by a result of computer experiments for $N=5$:

\begin{figure}[!htb]
\minipage{0.5\textwidth}
  \includegraphics[width=\linewidth]{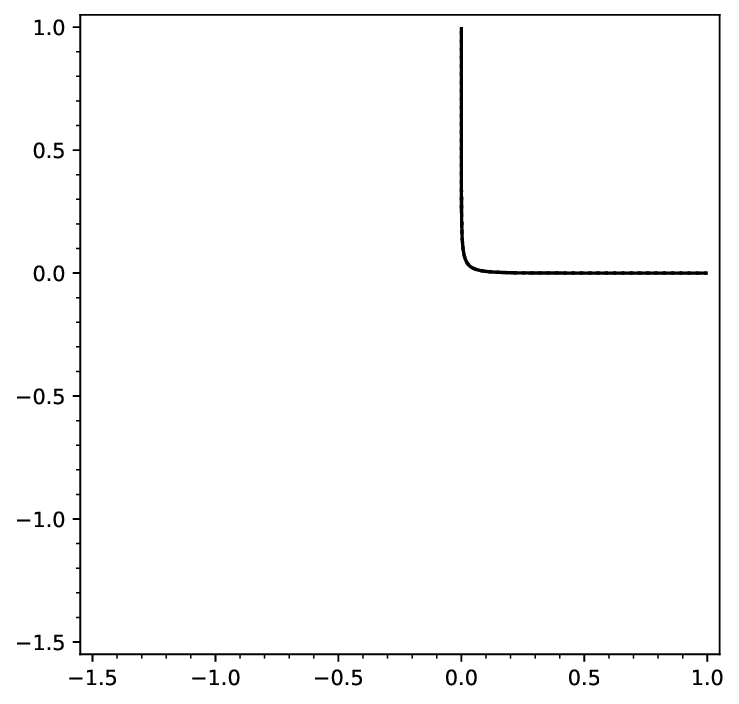}
  \caption*{$\epsilon=0$}
\endminipage\hfill
\minipage{0.5\textwidth}
  \includegraphics[width=\linewidth]{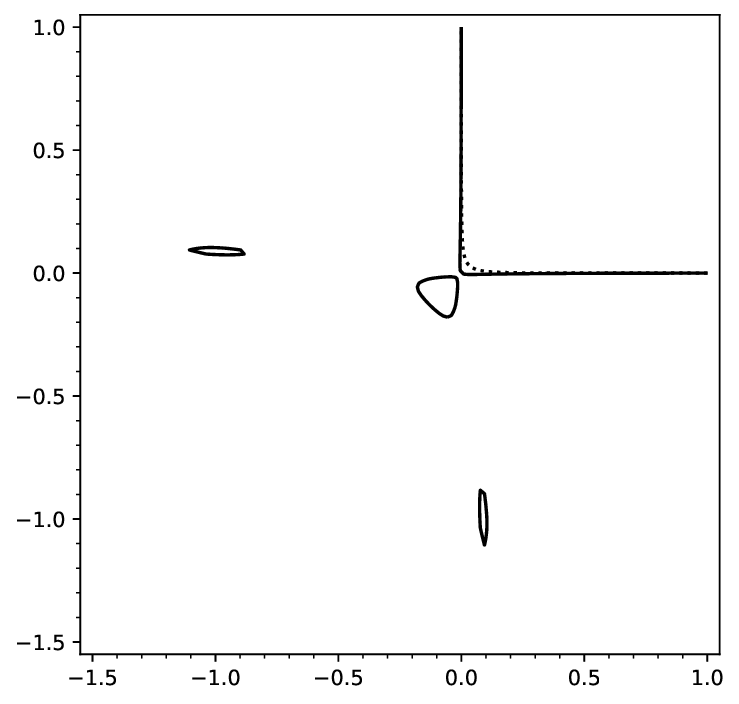}
  \caption*{$\epsilon=1$}
\endminipage\hfill
\vspace{0.5cm}
\minipage{0.5\textwidth}
  \includegraphics[width=\linewidth]{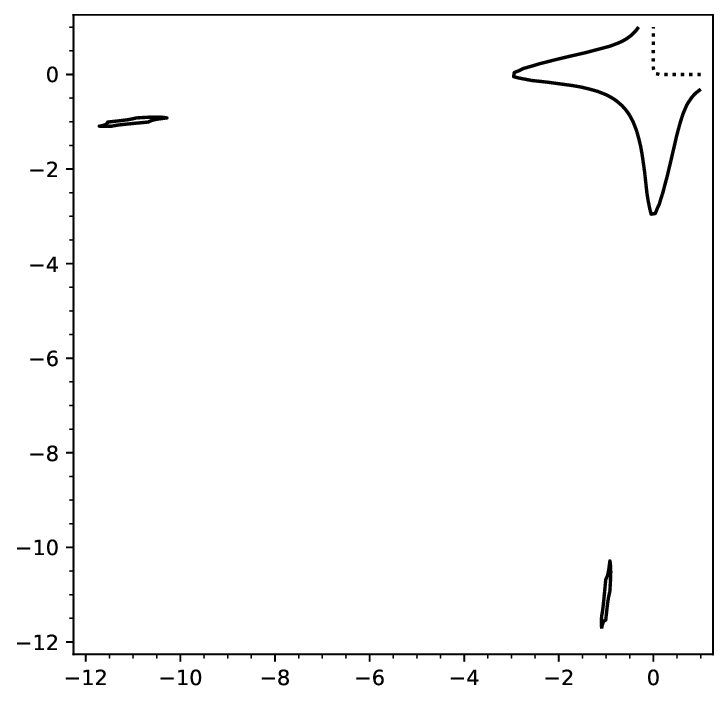}
  \caption*{$\epsilon=1000$}
  \endminipage\hfill
  \minipage{0.5\textwidth}
  \includegraphics[width=\linewidth]{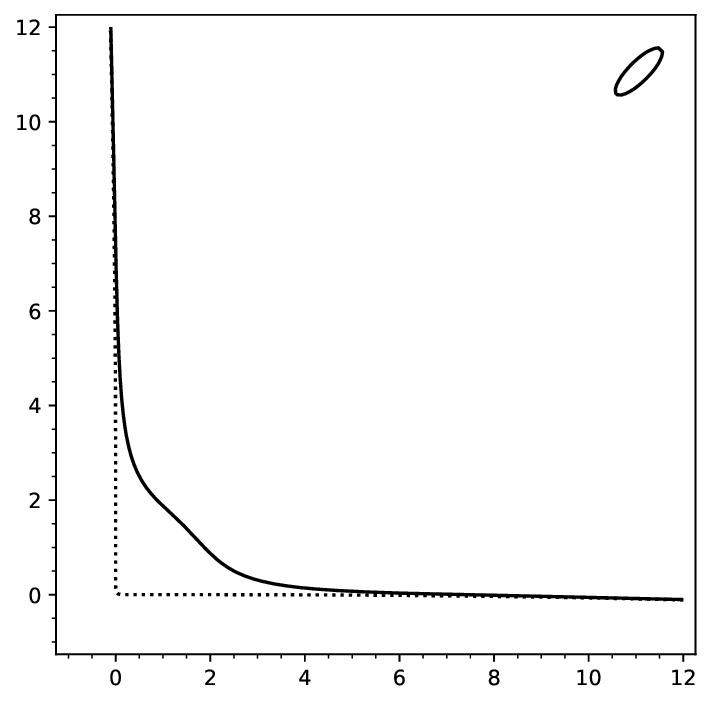}
  \caption*{$\epsilon=-1000$}
  \endminipage\hfill
\end{figure}

The ovals that appear shrink to the isolated  points when $\epsilon$ approaches to $0$.
\end{document}